\documentclass[11pt]{amsart}
\usepackage[latin1]{inputenc}
\usepackage[T1]{fontenc}
\usepackage{amsfonts}
\usepackage{amsmath}
\usepackage{amssymb}
\usepackage{graphicx}
\usepackage{pstricks}
\usepackage{pst-grad}
\usepackage{multido}
\usepackage{amsthm}
\usepackage{hyperref}
\usepackage{geometry}
\usepackage[absolute]{textpos}
\definecolor{purple}{rgb}{0.8,0.12,0.8}
\definecolor{orange}{rgb}{1.0,0.7,0.0}
\definecolor{pink}{rgb}{1,0.5,0.8}
\definecolor{blackg}{rgb}{0.1,0.25,0.1}
\definecolor{ForestGreen}{cmyk}{0.91,0,0.88,0.42}
\definecolor{Turquoise}{cmyk}{0.85,0,0.20,0}

\definecolor{GreenYellow}{cmyk}{0.15,0,0.69,0} 
\definecolor{Yellow}{cmyk}{0,0,1.,0} 
\definecolor{Goldenrod}{cmyk}{0,0.10,0.84,0} 
\definecolor{Dandelion}{cmyk}{0,0.29,0.84,0} 
\definecolor{Apricot}{cmyk}{0,0.32,0.52,0} 
\definecolor{Peach}{cmyk}{0,0.50,0.70,0} 
\definecolor{Melon}{cmyk}{0,0.46,0.50,0} 
\definecolor{YellowOrange}{cmyk}{0,0.42,1.,0} 
\definecolor{Orange}{cmyk}{0,0.61,0.87,0} 
\definecolor{BurntOrange}{cmyk}{0,0.51,1.,0} 
\definecolor{Bittersweet}{cmyk}{0,0.75,1.,0.24} 
\definecolor{RedOrange}{cmyk}{0,0.77,0.87,0} 
\definecolor{Mahogany}{cmyk}{0,0.85,0.87,0.35} 
\definecolor{Maroon}{cmyk}{0,0.87,0.68,0.32} 
\definecolor{BrickRed}{cmyk}{0,0.89,0.94,0.28} 
\definecolor{Red}{cmyk}{0,1.,1.,0} 
\definecolor{OrangeRed}{cmyk}{0,1.,0.50,0} 
\definecolor{RubineRed}{cmyk}{0,1.,0.13,0} 
\definecolor{WildStrawberry}{cmyk}{0,0.96,0.39,0} 
\definecolor{Salmon}{cmyk}{0,0.53,0.38,0} 
\definecolor{CarnationPink}{cmyk}{0,0.63,0,0} 
\definecolor{Magenta}{cmyk}{0,1.,0,0} 
\definecolor{VioletRed}{cmyk}{0,0.81,0,0} 
\definecolor{Rhodamine}{cmyk}{0,0.82,0,0} 
\definecolor{Mulberry}{cmyk}{0.34,0.90,0,0.02} 
\definecolor{RedViolet}{cmyk}{0.07,0.90,0,0.34} 
\definecolor{Fuchsia}{cmyk}{0.47,0.91,0,0.08} 
\definecolor{Lavender}{cmyk}{0,0.48,0,0} 
\definecolor{Thistle}{cmyk}{0.12,0.59,0,0} 
\definecolor{Orchid}{cmyk}{0.32,0.64,0,0} 
\definecolor{DarkOrchid}{cmyk}{0.40,0.80,0.20,0} 
\definecolor{Purple}{cmyk}{0.45,0.86,0,0} 
\definecolor{Plum}{cmyk}{0.50,1.,0,0} 
\definecolor{Violet}{cmyk}{0.79,0.88,0,0} 
\definecolor{RoyalPurple}{cmyk}{0.75,0.90,0,0} 
\definecolor{BlueViolet}{cmyk}{0.86,0.91,0,0.04} 
\definecolor{Periwinkle}{cmyk}{0.57,0.55,0,0} 
\definecolor{CadetBlue}{cmyk}{0.62,0.57,0.23,0} 
\definecolor{CornflowerBlue}{cmyk}{0.65,0.13,0,0} 
\definecolor{MidnightBlue}{cmyk}{0.98,0.13,0,0.43} 
\definecolor{NavyBlue}{cmyk}{0.94,0.54,0,0} 
\definecolor{RoyalBlue}{cmyk}{1.,0.50,0,0} 
\definecolor{Blue}{cmyk}{1.,1.,0,0} 
\definecolor{Cerulean}{cmyk}{0.94,0.11,0,0} 
\definecolor{Cyan}{cmyk}{1.,0,0,0} 
\definecolor{ProcessBlue}{cmyk}{0.96,0,0,0} 
\definecolor{SkyBlue}{cmyk}{0.62,0,0.12,0} 
\definecolor{Turquoise}{cmyk}{0.85,0,0.20,0} 
\definecolor{TealBlue}{cmyk}{0.86,0,0.34,0.02} 
\definecolor{Aquamarine}{cmyk}{0.82,0,0.30,0} 
\definecolor{BlueGreen}{cmyk}{0.85,0,0.33,0} 
\definecolor{Emerald}{cmyk}{1.,0,0.50,0} 
\definecolor{JungleGreen}{cmyk}{0.99,0,0.52,0} 
\definecolor{SeaGreen}{cmyk}{0.69,0,0.50,0} 
\definecolor{Green}{cmyk}{1.,0,1.,0} 
\definecolor{ForestGreen}{cmyk}{0.91,0,0.88,0.12} 
\definecolor{PineGreen}{cmyk}{0.92,0,0.59,0.25} 
\definecolor{LimeGreen}{cmyk}{0.50,0,1.,0} 
\definecolor{YellowGreen}{cmyk}{0.44,0,0.74,0} 
\definecolor{SpringGreen}{cmyk}{0.26,0,0.76,0} 
\definecolor{OliveGreen}{cmyk}{0.64,0,0.95,0.40} 
\definecolor{RawSienna }{cmyk}{0,0.72,1.,0.45} 
\definecolor{Sepia}{cmyk}{0,0.83,1.,0.70} 
\definecolor{Brown}{cmyk}{0,0.81,1.,0.60} 
\definecolor{Tan}{cmyk}{0.14,0.42,0.56,0} 
\definecolor{Gray}{cmyk}{0,0,0,0.50} 
\definecolor{Black}{cmyk}{0,0,0,1.} 
\definecolor{White}{cmyk}{0,0,0,0} 

\newcommand{\cA}{\mathcal{A}}
\newcommand{\cB}{\mathcal{B}}
\newcommand{\cC}{\mathcal{C}}
\newcommand{\cD}{\mathcal{D}}
\newcommand{\cE}{\mathcal{E}}
\newcommand{\cF}{\mathcal{F}}
\newcommand{\cG}{\mathcal{G}}
\newcommand{\cH}{\mathcal{H}}
\newcommand{\cI}{\mathcal{I}}
\newcommand{\cJ}{\mathcal{J}}
\newcommand{\cK}{\mathcal{K}}
\newcommand{\cL}{\mathcal{L}}
\newcommand{\cLR}{\mathcal{LR}}
\newcommand{\cM}{\mathcal{M}}
\newcommand{\cN}{\mathcal{N}}
\newcommand{\cO}{\mathcal{O}}
\newcommand{\cP}{\mathcal{P}}
\newcommand{\cQ}{\mathcal{Q}}
\newcommand{\cR}{\mathcal{R}}
\newcommand{\cS}{\mathcal{S}}
\newcommand{\cT}{\mathcal{T}}
\newcommand{\cU}{\mathcal{U}}
\newcommand{\cV}{\mathcal{V}}
\newcommand{\cW}{\mathcal{W}}
\newcommand{\cX}{\mathcal{X}}
\newcommand{\cY}{\mathcal{Y}}
\newcommand{\cZ}{\mathcal{Z}}

\newcommand{\sg}{\mathcal{h}}
\newcommand{\sd}{\mathcal{i}}

\newcommand{\bA}{\mathbf{A}}
\newcommand{\bB}{\mathbf{B}}
\newcommand{\bC}{\mathbf{C}}
\newcommand{\bD}{\mathbf{D}}
\newcommand{\bE}{\mathbf{E}}
\newcommand{\bF}{\mathbf{F}}
\newcommand{\bG}{\mathbf{G}}
\newcommand{\bH}{\mathbf{H}}
\newcommand{\bI}{\mathbf{I}}
\newcommand{\bJ}{\mathbf{J}}
\newcommand{\bK}{\mathbf{K}}
\newcommand{\bL}{\mathbf{L}}
\newcommand{\bM}{\mathbf{M}}
\newcommand{\bN}{\mathbf{N}}
\newcommand{\bO}{\mathbf{O}}
\newcommand{\bP}{\mathbf{P}}
\newcommand{\bQ}{\mathbf{Q}}
\newcommand{\bR}{\mathbf{R}}
\newcommand{\bS}{\mathbf{S}}
\newcommand{\bT}{\mathbf{T}}
\newcommand{\bU}{\mathbf{U}}
\newcommand{\bV}{\mathbf{V}}
\newcommand{\bW}{\mathbf{W}}
\newcommand{\bX}{\mathbf{X}}
\newcommand{\bY}{\mathbf{Y}}
\newcommand{\bZ}{\mathbf{Z}}

\newcommand{\ba}{\mathbf{a}}
\newcommand{\bb}{\mathbf{b}}
\newcommand{\bc}{\mathbf{c}}
\newcommand{\bd}{\mathbf{d}}
\newcommand{\bff}{\mathbf{f}}
\newcommand{\bg}{\mathbf{g}}
\newcommand{\bh}{\mathbf{h}}
\newcommand{\bj}{\mathbf{j}}
\newcommand{\bk}{\mathbf{k}}
\newcommand{\bl}{\mathbf{l}}
\newcommand{\bm}{\mathbf{m}}
\newcommand{\bn}{\mathbf{n}}
\newcommand{\bo}{\mathbf{o}}
\newcommand{\bp}{\mathbf{p}}
\newcommand{\bq}{\mathbf{q}}
\newcommand{\br}{\mathbf{r}}
\newcommand{\bs}{\mathbf{s}}
\newcommand{\bt}{\mathbf{t}}
\newcommand{\bu}{\mathbf{u}}
\newcommand{\bv}{\mathbf{v}}
\newcommand{\bw}{\mathbf{w}}
\newcommand{\bx}{\mathbf{x}}
\newcommand{\by}{\mathbf{y}}
\newcommand{\bz}{\mathbf{z}}

\newcommand{\fA}{\mathfrak{A}}
\newcommand{\fB}{\mathfrak{B}}
\newcommand{\fC}{\mathfrak{C}}
\newcommand{\fD}{\mathfrak{D}}
\newcommand{\fE}{\mathfrak{E}}
\newcommand{\fF}{\mathfrak{F}}
\newcommand{\fG}{\mathfrak{G}}
\newcommand{\fH}{\mathfrak{H}}
\newcommand{\fI}{\mathfrak{I}}
\newcommand{\fJ}{\mathfrak{J}}
\newcommand{\fK}{\mathfrak{K}}
\newcommand{\fL}{\mathfrak{L}}
\newcommand{\fM}{\mathfrak{M}}
\newcommand{\fN}{\mathfrak{N}}
\newcommand{\fO}{\mathfrak{O}}
\newcommand{\fP}{\mathfrak{P}}
\newcommand{\fQ}{\mathfrak{Q}}
\newcommand{\fR}{\mathfrak{R}}
\newcommand{\fS}{\mathfrak{S}}
\newcommand{\fT}{\mathfrak{T}}
\newcommand{\fU}{\mathfrak{U}}
\newcommand{\fV}{\mathfrak{V}}
\newcommand{\fW}{\mathfrak{W}}
\newcommand{\fX}{\mathfrak{X}}
\newcommand{\fY}{\mathfrak{Y}}
\newcommand{\fZ}{\mathfrak{Z}}

\newcommand{\nZ}{\mathbb{Z}}
\newcommand{\nR}{\mathbb{R}}
\newcommand{\nN}{\mathbb{N}}
\newcommand{\nQ}{\mathbb{Q}}
\newcommand{\nC}{\mathbb{C}}

\newcommand{\al}{\alpha}
\newcommand{\be}{\beta}
\newcommand{\si}{\sigma}
\newcommand{\la}{\lambda}
\newcommand{\ga}{\gamma}
\newcommand{\eps}{\epsilon}

\newcommand{\Ga}{\Gamma}
\newcommand{\De}{\Delta}
\newcommand{\Om}{\Omega}

\newcommand{\tA}{\tilde{A}}
\newcommand{\tB}{\tilde{B}}
\newcommand{\tC}{\tilde{C}}
\newcommand{\tD}{\tilde{D}}
\newcommand{\tE}{\tilde{E}}
\newcommand{\tF}{\tilde{F}}
\newcommand{\tG}{\tilde{G}}
\newcommand{\tW}{\tilde{W}}
\newcommand{\tJ}{\tilde{J}}
\newcommand{\tL}{\tilde{L}}
\newcommand{\tcH}{\tilde{\cH}}
\newcommand{\tc}{\tilde{c}}
\newcommand{\tb}{\tilde{b}}
\newcommand{\tT}{\tilde{t}}

\newcommand{\sqs}{\sqsubset}
\newcommand{\sq}{\sqsubseteq}
\newcommand{\ov}{\overline}
\newcommand{\ind}{\underset}
\newcommand{\sD}{\overline{\mathcal{F}}}
\newcommand{\tx}{\quad\text}
\newcommand{\mand}{\quad\text{and}\quad}

\newcommand{\Hi}{H^{(i)}}
\newcommand{\Hj}{H^{(j)}}

\newcommand{\ra}{\rightarrow}
\newcommand{\lra}{\longrightarrow}
\newcommand{\Ra}{\Rightarrow}
\newcommand{\lar}{\leftarrow}
\newcommand{\llar}{\longleftarrow}
\newcommand{\Lar}{\Leftarrow}
\newcommand{\eq}{\Leftrightarrow}

\newenvironment{longlist}%
{ \begin{list}%
	{$\bullet$}%
	{\setlength{\labelwidth}{30pt}%
	 \setlength{\leftmargin}{0pt}%
	 \setlength{\itemsep}{.1cm}}}%
{ \end{list} }

\newtheorem{Th}{Theorem}[section]

\newtheorem{Def-Prop}[Th]{Definition-Proposition}

\theoremstyle{definition}

\newtheorem{Cl}[Th]{Claim}

\theoremstyle{remark}
\newtheorem{Rem}[Th]{Remark}

\begin{document}
\bibliographystyle{plain}
\title{Some computations about Kazhdan-Lusztig cells in affine Weyl groups of rank 2}
\author{J\'er\'emie Guilhot}

\address{School of Mathematics and Statistics F07
University of Sydney NSW 2006
Australia }

\email{guilhot@maths.usyd.edu.au}


\maketitle
\begin{abstract}
In the last section of the paper ``Generalized induction of Kazhdan-Lusztig cells'' and in ``Kazhdan-Lusztig cells in affine Weyl groups of rank 2''  the author described the partition into Kazhdan-Lusztig cells of the affine Weyl groups of rank 2 for all choices of parameters. The proof of these results relies on some explicit computation with GAP \cite{gap}. In these notes we give some detail of these computations. \end{abstract}

\section{Introduction}

These notes are organized as follows:

\begin{itemize}
\item in Section \ref{sec1}, we give some details in the computations involved in the last section of \cite[\S 6]{jeju3};
\item in Section \ref{sec2}, we give some details in the computations involved in the proof of \cite[Theorem 5.1]{jeju4} for the affine Weyl group of type $\tG_{2}$. We also give the left order on the left cells for all choices of parameters.
\item in Section \ref{sec3}, we describe the partition of  the affine Weyl group of type $\tB_{2}$ into cells for all parameters.
\item in Section \ref{sec4}, we give some details in the computations involved in the proof of \cite[Theorem 5.3]{jeju4} for the affine Weyl group of type $\tB_{2}$. We also give the left order on the left cells for all choices of parameters.

\end{itemize}

\section{Computations in \cite[\S 6]{jeju3}}
\label{sec1}
Let $W$ be the affine Weyl group of type $\tilde{G_{2}}$ with diagram and weight function given by
\begin{center}

\end{center}
where $a,b$ are positive integers such that $a/b>4$. In this section we go through the proof of Lemma 6.4, Claim 6.6 and Claim 6.8 in \cite{jeju3} and give some details in the computations involved. \\
$\ $\\

We keep the notation of \cite[\S 6]{jeju3}. Recall that
$\ $\\
\begin{center}
%
\end{center}

\subsection{The sets $C_{i}$} $\ $\\

For $1\leq i\leq 6$, let
\begin{enumerate}
\item $u_{i}\in C_{i}$ be the element of minimal length in $C_{i}$;
\item $v_{i}\in A_{i}$ be the element of minimal length in $A_{i}$;
\item $v'_{i}\in A'_{i}$ be the element of minimal length in $A'_{i}$.
\end{enumerate}
We set $U:=\{u_{i},v_{i},v'_{i}\ |\ 1\leq i\leq 6\}$, $X_{v_{i}}=X_{v'_{i}}=X_{\{s_{1},s_{2}\}}$ and
$$X_{u_{i}}=\{z\in W\ |\ zu_{i}\in C_{i}\}.$$
In this section we prove the following.
\begin{enumerate}
\item The submodule
$$\cM:=\sg T_{x}C_{u}\ |\ u\in U, x\in X_{u}\sd_{\cA}\subseteq \cH$$
is a left ideal.
\item Let $u\in U$, $u'<u$ and $x\in X_{u}$. We have
$$P_{u',u}T_{x}C_{u}\equiv T_{xu} \mod \cH_{<0}$$
\end{enumerate}


\noindent First we prove (1). To this end we only need to show that $T_{s_{2}}C_{u_{i}}$ and $T_{s_{2}s_{1}s_{2}s_{1}s_{2}s_{3}}C_{u_{i}}$ lie in $\cM$ for all $1\leq i\leq 6$. We have
$$T_{s_{2}}C_{u_{i}}=C_{v_{i}}-v^{-b}C_{u_{i}}\tx{for $i=1,2,3,6$}$$
and
$$T_{s_{2}}C_{u_{4}}=C_{v_{4}}-v^{-b}C_{u_{4}}+C_{v_{2}}$$
$$T_{s_{2}}C_{u_{5}}=C_{v_{5}}-v^{-b}C_{u_{5}}+C_{v_{1}}.$$
Thus $T_{s_{2}}C_{u_{i}}\in\cM$.\\
For $i=1,2,3,6$ we have
\begin{eqnarray*}
T_{3}C_{u_{i}}&=&C_{3u_{i}}+\cA C_{u_{i}}\\
T_{23}C_{u_{i}}&=&C_{23u_{i}}+\cA C_{3u_{i}}+\cA C_{u_{i}}+\cA C_{v_{i}}\\
T_{123}C_{u_{i}}&=&C_{123u_{i}}+\cA C_{23u_{i}}+\cA C_{3u_{i}}+\cA C_{u_{i}}+\cA C_{v_{i}}\\
T_{2123}C_{u_{i}}&=&C_{2123u_{i}}+\cA C_{123u_{i}}+\cA C_{3u_{i}}+\cA C_{u_{i}}\\
T_{12123}C_{u_{i}}&=&C_{12123u_{i}}+\cA C_{2123u_{i}}+\cA C_{123u_{i}}+\cA C_{u_{i}}\\
T_{212123}C_{u_{i}}&=&C_{v'_{i}}+\cA C_{12123u_{i}}+\cA C_{123u_{i}}+\cA C_{u_{i}}+\cA C_{v_{i}}\\
\end{eqnarray*}
thus 
$$C_{3u_{i}},C_{23u_{i}},C_{123u_{i}},C_{2123u_{i}},C_{12123u_{i}}\in\cM.$$
Since $C_{v'_{i}},C_{u_{i}}, C_{v_{i}}\in\cM$ we see that $T_{212123}C_{u_{i}}\in\cM$ as required.\\
We have
\begin{eqnarray*}
T_{3}C_{u_{4}}&=&C_{3u_{4}}+\cA C_{u_{4}}\\
T_{23}C_{u_{4}}&=&C_{23u_{4}}+\cA C_{3u_{4}}+\cA C_{u_{4}}+\cA C_{v_{4}}+\cA C_{v_{2}}\\
T_{123}C_{u_{4}}&=&C_{123u_{4}}+\cA C_{23u_{4}}+\cA C_{3u_{4}}+\cA C_{u_{4}}+\cA C_{v_{4}}+\cA C_{v_{2}}\\
T_{2123}C_{u_{4}}&=&C_{2123u_{4}}+\cA C_{123u_{4}}+\cA C_{3u_{4}}+\cA C_{u_{4}}\\
T_{12123}C_{u_{4}}&=&C_{12123u_{4}}+\cA C_{2123u_{4}}+\cA C_{123u_{4}}+\cA C_{u_{4}}\\
T_{212123}C_{u_{4}}&=&C_{v'_{i}}+\cA C_{12123u_{4}}+\cA C_{123u_{4}}+\cA C_{u_{4}}+\cA C_{v_{4}}+\cA C_{v_{2}}\\
\end{eqnarray*}
thus 
$$C_{3u_{i}},C_{23u_{i}},C_{123u_{i}},C_{2123u_{i}},C_{12123u_{i}}\in\cM.$$
Since $C_{v'_{i}},C_{u_{i}}, C_{v_{i}}\in\cM$ we see that $T_{212123}C_{u_{4}}\in\cM$ as required.\\
We have
\begin{eqnarray*}
T_{3}C_{u_{5}}&=&C_{3u_{5}}+\cA C_{u_{5}}\\
T_{23}C_{u_{5}}&=&C_{23u_{5}}+\cA C_{3u_{5}}+\cA C_{u_{5}}+\cA C_{v_{5}}+\cA C_{v_{1}}\\
T_{123}C_{u_{5}}&=&C_{123u_{5}}+\cA C_{23u_{5}}+\cA C_{3u_{5}}+\cA C_{u_{5}}+\cA C_{v_{5}}+\cA C_{v_{1}}\\
T_{2123}C_{u_{5}}&=&C_{2123u_{5}}+\cA C_{123u_{5}}+\cA C_{3u_{5}}+\cA C_{u_{5}}\\
T_{12123}C_{u_{5}}&=&C_{12123u_{5}}+\cA C_{2123u_{5}}+\cA C_{123u_{5}}+\cA C_{u_{5}}\\
T_{212123}C_{u_{5}}&=&C_{v'_{i}}+\cA C_{12123u_{5}}+\cA C_{123u_{5}}+\cA C_{u_{5}}+\cA C_{v_{5}}+\cA C_{v_{1}}\\
\end{eqnarray*}
thus 
$$C_{3u_{i}},C_{23u_{i}},C_{123u_{i}},C_{2123u_{i}},C_{12123u_{i}}\in\cM.$$
Since $C_{v'_{i}},C_{u_{i}}, C_{v_{i}}\in\cM$ we see that $T_{212123}C_{u_{5}}\in\cM$ as required.\\
Statement (1) follows.

\vspace{.5cm}
\noindent
Before proving (2), we need to introduce some definitions. For $w\in W$ we set
$$Y(w)=\{y\in W | \ell(wy^{-1})=\ell(w)-\ell(y)\}.$$
We have $y\in Y(w)$ if and only if there exists a reduced expression $s_{1}\ldots s_{n}$ of $w$ and $1\leq i\leq n$ such that $y=s_{i}\ldots s_{n}$. Recall the definition of translation in \cite[Definition 4.1]{jeju1}.
\begin{Cl}
\label{trans}
Let $w\in W$ be a translation. Let $z\in W$. The following statements are equivalent
\begin{enumerate}
\item $\ell(wz)=\ell(w)+\ell(z)$;
\item $\ell(xz)=\ell(x)+\ell(z) \text{ for all $x\in \{zw^{n}| z\in Y(w), n\in\nN \}$}$.
\end{enumerate}
\end{Cl}
\begin{proof}
The fact that (2) implies (1) is clear since $w\in \{zw^{n}| z\in Y(w), n\in\nN \}$. Assume that (1) holds. Then using \cite[Lemma 4.2]{jeju1} in its right version we get
$$\ell(wz)=\ell(w)+\ell(z)\Leftrightarrow \ell(w^{n}z)=\ell(w^{n})+\ell(z).$$
Let $x\in  \{zw^{n}| z\in Y(w), n\in\nN \}$. Let $y\in Y_{w}$ and $n\in \nN$ be such that $x=yw^{n}$. We have
$$\ell(wz)=\ell(w)+\ell(z)\Rightarrow \ell(w^{n+1}z)=\ell(w^{n+1})+\ell(z)\Rightarrow \ell(yw^{n}z)=\ell(yw^{n})+\ell(z)$$
as required.
\end{proof}

\vspace{.5cm}
We now prove (2).  Let $w=s_{1}s_{2}s_{1}s_{2}s_{3}s_{1}s_{2}s_{1}s_{2}s_{3}$. Note that $w$ is a translation. We have
$$X_{u_{i}}=\{yw^{n}| y\in Y(w), n\in\nN\}.$$
Let $u<u_{i}$. Using GAP one can show that for all $y\in Y(w)$ we have 
$$P_{u,u_{i}}T_{y}C_{u}\equiv T_{yu} \mod \cH_{<0}.$$
Next, one can show that 
$$P_{u,u_{i}}T_{w}C_{u}=T_{wu}+\sum_{z\in W} a_{z}T_{z}$$
where $a_{z}\in \cA_{<0}$ and if $a_{z}\neq 0$ then $z$ satisfy $\ell(wz)=\ell(w)+\ell(z)$. Now let $x=yw^{n}\in X_{u_{i}}$ (where $y\in Y(w)$ and $n\in\nN^{*}$). Then using Claim \ref{trans} we obtain
$$\begin{array}{ccl}
P_{u,u_{i}}T_{x}C_{u}&=&T_{yw^{n-1}}\bigg(T_{wu}+\sum_{z\in W} a_{z}T_{z}\bigg)\\
&=& T_{xu}+\underset{z\in W}{\sum}  a_{z}T_{yw^{n-1}z}\\
&\equiv& T_{xu} \mod \cH_{<0}
\end{array}$$
and statement (2) follows. 
$\ $\\

\subsection{Proof of Claim 6.6} $\ $\\

\noindent Let $u=s_{1}s_{3}s_{2}s_{1}$ and recall that 
\begin{align*}
u_{1}&=s_{1}s_{2}s_{1}s_{2}s_{1},\\
v_{1}&=s_{1}s_{2}s_{1}s_{2}s_{1}s_{2},\\
v'_{1}&=s_{1}s_{2}s_{1}s_{2}s_{1}s_{2}s_{3}s_{2}s_{1}s_{2}s_{1},\\
u_{2}&=s_{1}s_{2}s_{1}s_{2}s_{1}s_{3}s_{2}s_{1},\\
v_{2}&=s_{2}s_{1}s_{2}s_{1}s_{2}s_{1}s_{3}s_{2}s_{1},\\
v'_{2}&=s_{2}s_{1}s_{2}s_{1}s_{2}s_{3}s_{1}s_{2}s_{1}s_{2}s_{1}s_{3}s_{2}s_{1}\\
v_{3}&=s_{2}s_{1}s_{2}s_{1}s_{2}s_{1}s_{3},
\end{align*}
and
\begin{align*}
X_{u}&:=\{z\in W| zu\in B_{1}\}\\
X_{u_{i}}&=\{z\in W\ |\ zu_{i}\in C_{i}\},\\
X_{v_{i}}&=X_{v'_{i}}=X_{\{s_{1},s_{2}\}}\quad\text{for $1\leq i\leq 6$}.
\end{align*}
Let $U:=\{u,u_{1},v_{1},v'_{1},u_{2},v_{2},v'_{2},v_{3}\}$. In this section we prove the following.
\begin{enumerate}
\item The submodule
$$\cM:=\sg T_{x}C_{v}\ |\ v\in U, x\in X_{u}\sd_{\cA}\subseteq \cH$$
is a left ideal.
\item Let $u'<u$ and $x\in X_{u}$. We have
$$P_{u',u}T_{x}C_{u}\equiv T_{xu} \mod \cH_{<0}$$
\end{enumerate}

\vspace{.5cm}
\noindent We prove (1). To this end it is enough to show that $T_{1212}C_{u}\in\cM$ and $T_{1212132}C_{u}$ lie in $\cM$. We have
\begin{eqnarray*}
T_{2}C_{u}&=&C_{2u}+\cA C_{u}\\
T_{12}C_{u}&=&C_{12u}+ \cA C_{2u}+\cA C_{u}\\
T_{212}C_{u}&=&C_{212u}+\cA C_{12u}+\cA C_{u}\\
T_{1212}C_{u}&=&C_{u_{2}}+ \cA C_{212u}+\cA C_{12u}+\cA C_{u_{1}}+\cA C_{v_{1}}
\end{eqnarray*}
Thus $T_{1212}C_{u}\in \cM$. \\
\newpage
We have
\begin{eqnarray*}
T_{2}C_{u}&=&C_{2u}+\cA C_{u}\\
T_{32}C_{u}&=&C_{32u}+\cA C_{2u}\\
T_{132}C_{u}&=&C_{132u}+ \cA C_{32u}+\cA C_{12u}+\cA C_{2u}+ \cA C_{u}\\
T_{2132}C_{u}&=&C_{2132u}+ \cA C_{212u}+\cA C_{132u}+\cA C_{32u}+\cA C_{12u}+\cA C_{2u}+ \cA C_{u}\\
T_{12132}C_{u}&=&C_{12132u}+ \cA C_{2132u}+\cA C_{212u}+\cA C_{32u}+\cA C_{12u}+\cA C_{2u}\\
&&+ \cA C_{u}+\cA C_{u_{1}}+\cA C_{v_{1}}+\cA C_{u_{2}}\\
T_{212132}C_{u}&=&C_{212132u}+\cA C_{12132u}+ \cA C_{2132u}+\cA C_{212u}+\cA C_{32u}\\ &&+\cA C_{12u}+\cA C_{2u}+ \cA C_{u}+\cA C_{u_{1}}+\cA C_{v_{2}}+\cA C_{u_{2}}\\
T_{1212132}C_{u}&=&C_{v'_{1}}+\cA C_{212132u}+\cA C_{2132u} +\cA C_{212u}+\cA C_{32u}+\cA C_{2u}\\
&& +\cA C_{v_{1}}+\cA C_{v_{2}}+\cA C_{v_{3}}
\end{eqnarray*}
Thus $T_{1212132}C_{u}\in\cM$ as desired.

\vspace{1cm}

\noindent We prove (2). Let $w_{1}=s_{1}s_{2}s_{3}s_{2}s_{1}s_{2}$ and $w_{2}=s_{2}s_{3}s_{2}s_{1}s_{2}s_{1}$. These are both translations. Any $x\in X_{u}$ can be written under the form $yw_{1}^{n}$ with $y\in Y(w_{1})$ or $yw_{2}^{n}s_{3}s_{2}$ with $y\in Y(w_{2})$. One can show that for all $y_{1}\in Y(w_{1})$ and all $y_{2}\in Y(w_{2})$ we have
$$P_{u,u_{i}}T_{y_{1}}C_{u}\equiv T_{y_{1}u} \mod \cH_{<0}\text{ and }  P_{u,u_{i}}T_{y_{2}s_{3}s_{2}}C_{u}\equiv T_{y_{2}s_{3}s_{2}u} \mod \cH_{<0}$$
Furthermore, we have
$$P_{u,u_{i}}T_{w_{1}}C_{u}=T_{w_{1}u}+\sum_{z\in W} a_{z}T_{z}$$
where $a_{z}\in\cA_{<0}$ and if $a_{z}\neq 0$ then $z$ satisfy $\ell(w_{1}z)=\ell(w_{1})+\ell(z)$. We have 
$$P_{u,u_{i}}T_{w_{2}s_{3}s_{2}}C_{u}=T_{w_{2}s_{3}s_{2}u}+\sum_{z\in W} a'_{z}T_{z}$$
where $a'_{z}\in\cA_{<0}$ and if $a'_{z}\neq 0$ then $z$ satisfy $\ell(w_{2}z)=\ell(w_{2})+\ell(z)$. 
Arguing as in the previous section we obtain that (2) holds.\\

\subsection{Proof of Claim 6.8} $\ $\\

\noindent Recall that $v=s_{1}s_{3}s_{2}s_{1}s_{2}s_{3}$ and 
$$X_{v}:=\{z\in W|zv\in B_{3}\}\qquad Y_{v}:=\{x\in X_{v}| \ell(xs_{2}s_{1}s_{2})=\ell(x)-3\}.$$
Let $U=\{v,u_{4},v_{4},v'_{4},v_{3},v_{2},v_{5}\}$.  In this section we prove the following.
\begin{enumerate}
\item The submodule
$$\cM:=\sg T_{x}C_{u}\ |\ u\in U, x\in X_{u}\sd_{\cA}\subseteq \cH$$
is a left ideal.
\item Let $x\in X_{v}-Y_{v}$. Then we have
$$T_{x}C_{v}\equiv T_{xv} \mod \cH_{<0}$$
\item Let $y\in Y_{v}$ and $y_{0}=s_{2}s_{1}s_{2}$. We have
$$\tC_{yv}\equiv T_{yv}+T_{y_{0}v_{3}} \mod \cH_{<0}$$

\end{enumerate}

\vspace{.5cm}

\noindent We prove (1). To this end it is enough to show that $T_{1212}C_{u}\in\cM$ and $T_{1212132}C_{u}$ lie in $\cM$. We have
\begin{eqnarray*}
T_{2}C_{v}&=&C_{2v}+\cA C_{v}\\
T_{12}C_{v}&=&C_{12v}+ \cA C_{2v}+\cA C_{v}\\
T_{212}C_{v}&=&C_{212v}+\cA C_{12v}+ \cA C_{v}+\cA C_{v_{3}}\\
T_{1212}C_{v}&=&C_{u_{4}}+ \cA C_{212v}+\cA C_{12v}+\cA C_{v_{3}}
\end{eqnarray*}
Thus $T_{1212}C_{v}\in \cM$. We have
\begin{eqnarray*}
T_{2}C_{v}&=&C_{2v}+\cA C_{v}\\
T_{32}C_{v}&=&C_{32v}+\cA C_{2v}\\
T_{132}C_{v}&=&C_{132v}+ \cA C_{32v}+\cA C_{12v}+\cA C_{2v}+ \cA C_{v}\\
T_{2132}C_{v}&=&C_{2132v}+ \cA C_{212v}+\cA C_{132v}+\cA C_{32v}+\cA C_{12v}+\cA C_{2v}+ \cA C_{v}+\cA C_{v_{3}}\\
T_{12132}C_{v}&=&C_{12132v}+ \cA C_{2132v}+\cA C_{212v}+\cA C_{32v}+\cA C_{12v}+\cA C_{2v}\\
&&+ \cA C_{v}+\cA C_{u_{4}}+\cA C_{v_{3}}\\
T_{212132}C_{v}&=&C_{212132v}+\cA C_{12132v}+ \cA C_{2132v}+\cA C_{212v}+\cA C_{32v}\\ &&+\cA C_{12v}+\cA C_{2v}+ \cA C_{v}+\cA C_{u_{4}}+\cA C_{v_{4}}+\cA C_{v_{3}}+\cA C_{v_{2}}\\
T_{1212132}C_{v}&=&T_{212123}C_{v_{3}}+\cA C_{212132v}+\cA C_{2132} +\cA C_{212v}+\cA C_{32v}+\cA C_{2v}\\
&& +\cA C_{v_{5}}+\cA C_{v_{4}}+\cA C_{v_{3}}+\cA C_{v_{2}}
\end{eqnarray*}
Thus $T_{1212132}C_{u}\in\cM$ as desired.\\

\vspace{1cm}

\noindent We prove (2) and (3).
Let $w_{1}=s_{1}s_{2}s_{3}s_{2}s_{1}s_{2}$ and $w_{2}=s_{2}s_{3}s_{2}s_{1}s_{2}s_{1}$. These are both translations. Any element in $X_{v}$ can be written under the form $y_{1}w_{1}^{n}$ or $y_{2}w_{2}^{n}$ (where $y_{1}\in Y(w_{1})$ and $y_{2}\in Y(w_{2})$). 
Now we have
$$\begin{array}{rllcc}
T_{s_{2}}C_{v}&\equiv& T_{s_{2}v} &\mod \cH_{<0}\\ 
T_{s_{1}s_{2}}C_{v}&\equiv& T_{s_{1}s_{2}v} & \mod \cH_{<0}\\ 
T_{s_{2}s_{1}s_{2}}C_{v}&\equiv& T_{s_{2}s_{1}s_{2}v}+T_{v_{3}}& \mod \cH_{<0}\\ 
T_{s_{3}s_{2}s_{1}s_{2}}C_{v}&\equiv& T_{s_{3}s_{2}s_{1}s_{2}v}+T_{s_{3}v_{3}}& \mod \cH_{<0}\\ 
T_{s_{2}s_{3}s_{2}s_{1}s_{2}}C_{v}&\equiv& T_{s_{2}s_{3}s_{2}s_{1}s_{2}v}+T_{s_{2}s_{3}v_{3}}&  \mod \cH_{<0}\\ 
T_{w_{1}}C_{v}&\equiv& T_{w_{1}v}+T_{s_{1}s_{2}s_{3}v_{3}}&  \mod \cH_{<0}\\ 
\end{array}$$
and if $T_{z}$ appears with a non zero coefficient in the expression of $T_{w_{1}}C_{v}$ in the standard basis then $\ell(w_{1}z)=\ell(w_{1})+\ell(z)$. Therefore for all $y_{1}\in Y(w_{1})$ and $n\geq 1$ we have
$$T_{y_{1}w_{1}^{n}}C_{v} \equiv T_{yw_{1}^{n-1}s_{2}s_{1}s_{2}v}+T_{y_{1}w_{1}^{n-1}s_{1}s_{2}s_{3}v_{3}}  \mod \cH_{<0}.$$
Next we have
$$\begin{array}{rllcc}
T_{s_{2}}C_{v}&\equiv& T_{s_{2}v} &\mod \cH_{<0}\\ 
T_{s_{3}s_{2}}C_{v}&\equiv& T_{s_{3}s_{2}v} & \mod \cH_{<0}\\ 
T_{s_{1}s_{3}s_{2}}C_{v}&\equiv& T_{s_{1}s_{3}s_{2}v}& \mod \cH_{<0}\\ 
T_{s_{2}s_{1}s_{3}s_{2}}C_{v}&\equiv&T_{s_{2}s_{1}s_{3}s_{2}v}& \mod \cH_{<0}\\ 
T_{s_{1}s_{2}s_{1}s_{3}s_{2}}C_{v}&\equiv& T_{s_{1}s_{2}s_{1}s_{3}s_{2}v} & \mod \cH_{<0}\\ 
T_{s_{2}s_{1}s_{2}s_{1}s_{3}s_{2}}C_{v}&\equiv& T_{s_{2}s_{1}s_{2}s_{1}s_{3}s_{2}v}& \mod \cH_{<0}\\ 
T_{s_{3}s_{2}s_{1}s_{2}s_{1}s_{3}s_{2}}C_{v}&\equiv&  T_{s_{3}s_{2}s_{1}s_{2}s_{1}s_{3}s_{2}v}& \mod \cH_{<0}\\ 
T_{w_{2}s_{3}s_{2}}C_{v}&\equiv& T_{w_{2}v}+T_{s_{3}s_{2}s_{1}s_{2}s_{3}v_{3}} & \mod \cH_{<0}\\ 
\end{array}$$
and if $T_{z}$ appears with a non zero coefficient in the expression of $T_{w_{2}s_{3}s_{2}}C_{v}$ in the standard basis then $\ell(w_{2}z)=\ell(w_{2})+\ell(z)$.  Therefore for all $y_{2}\in Y(w_{2})$ and $n\geq 1$ we have
$$T_{y_{2}w_{2}^{n}s_{3}s_{2}}C_{v} \equiv T_{y_{2}w_{2}^{n-1}s_{2}s_{1}s_{2}v}+T_{y_{2}w_{2}^{n-1}s_{3}s_{2}s_{1}s_{2}s_{3}v_{3}}  \mod \cH_{<0}.$$
Let $y\in Y_{v}$ and let $y_{0}=ys_{2}s_{1}s_{2}$. Using the previous equalities, we obtain that 
$$T_{y}C_{v}\equiv T_{yv}+T_{y_{0}v_{3}} \mod \cH_{<0}.$$
and for $x\in X_{v}-Y_{v}$ we have
$$T_{y}C_{v}\equiv T_{yv} \mod \cH_{<0}.$$
Finally using the fact that $p^{*}_{xu,yv}\in\cA_{<0}$ we get
$$
\begin{array}{rllll}
\tilde{C_{yv}}&=&T_{y}C_{v}+\ind{xu\sqs yv}{\ind{u\in U, x\in X_{u}}{\sum}}p^{*}_{xu,yv}T_{x}C_{u}&\\
&=&T_{y}C_{v}+\ind{x< y,x\in X_{v}}{\sum}p^{*}_{xv,yv}T_{x}C_{v}+ \ind{u\neq v}{\ind{u\in U, x\in X_{u}}{\sum}}p^{*}_{xu,yv}T_{x}C_{u}&\\
&=&T_{y}C_{v} \qquad\qquad\qquad\qquad\qquad\qquad\  &\text{mod $\cH_{<0}$}\\
&=&T_{yv}+T_{y_{0}v_{3}}& \text{mod $\cH_{<0}$}\\
\end{array}
$$
as required.

\newpage

\section{Proof of \cite[Theorem 5.1]{jeju4} for $\tG_{2}$}
\label{sec2}

We keep the notation of \cite{jeju4}.\\

\noindent
Let $(W,S)$ be the affine Weyl group of type $\tG_{2}$ with diagram and weight function given by
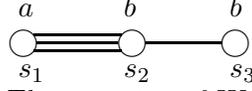
\begin{center}
\begin{picture}(150,32)
\put( 40, 10){\circle{10}}
\put( 44,  7){\line(1,0){33}}
\put( 45,  10){\line(1,0){30.5}}
\put( 44, 13){\line(1,0){33}}
\put( 81, 10){\circle{10}}
\put( 86, 10){\line(1,0){29}}
\put(120, 10){\circle{10}}
\put( 38, 20){$a$}
\put( 78, 20){$b$}
\put(118, 20){$b$}
\put( 38, -3){$s_{1}$}
\put( 78, -3){$s_{2}$}
\put(118,-3){$s_{3}$}
\end{picture}
\end{center}
where $a,b$ are positive integers. The partition of $W$ into cells only depends on the ratio $r=a/b$. The aim of this section is to give some details in the computation involved in \cite[Theorem 5.1]{jeju4}. Recall that in \cite{jeju4} we showed that the partition of $W$ given in \cite[Conjecture 3.3]{jeju4} is constant whether
\begin{itemize}
\item $r>2$
\item $r=2$
\item $2>r>3/2$
\item $r=3/2$
\item $3/2>r>1$
\item $r=1$
\item $r<1$
\end{itemize}

Fix $r>0$. Let $\tc_{0},\ldots,\tc_{n}$ the sets obtained using the algorithm presented in \cite[\S 3]{jeju4}. Then we denote by
\begin{itemize}
\item $\tc_{i}^{j}$ the left-connected components lying in $\tc_{i}$;
\item $u_{i}^{j}$ the element of minimal length lying in $\tc_{i}^{j}$;
\item $U$ the set which consists of all the $u_{i}^{j}$.
\end{itemize}
Let $u_{i}^{j}\in U$. We set
$$X_{u_{i}^{j}}:=\{w\in W| wu^{j}_{i}\in \tc^{j}_{i}\}.$$
For each value of $r$ we give
\begin{enumerate}
\item the Hasse diagram of the pre-order $\preceq$ on $U$, which is also the left order on the left cells;
\item A table where for $u\in U$ and $x\in X_{u}$ we put the expression 
$$T_{yv}+\sum_{xu\sqs yv} a_{xu,yv}T_{x}C_{u}$$
such that
$$T_{y}C_{v}\equiv T_{yv}+\sum_{xu\sqs yv} a_{xu,yv}T_{x}C_{u} \mod \cH_{<0} .$$
\begin{Rem}
When $T_{x}C_{u}\equiv T_{xu} \mod \cH_{<0}$ we do not put it in the table.
\end{Rem}
\end{enumerate}

\noindent
For $w\in W$ we set
$$Y(w)=\{y\in W| \ell(wy^{-1})=\ell(w)-\ell(y)\}.$$

\newpage

\begin{longlist}


\item {\bf Case ${\bf r>2}$.}\\

\begin{center}

\end{center}

\newpage 

In the following Hasse diagram, the black edges hold for all $r>2$, the green edges only hold for $3<r\leq 4$ and the red edges only hold for  $2<r\leq 3$.

\begin{center}
%
\end{center}

\newpage

Let $w=s_{1}s_{2}s_{1}s_{2}s_{3}s_{1}s_{2}s_{1}s_{2}s_{3}$, $w_{1}=s_{1}s_{2}s_{3}s_{2}s_{1}s_{2}$ and $w_{2}=s_{2}s_{3}s_{2}s_{1}s_{2}s_{1}$. In the following table we have
$$y_{1}\in \{yw_{1}^{n}| y\in Y(w_{1}), n\in\nN\} \text{ and } y_{2}\in \{yw_{2}^{n}| y\in Y(w_{2})\}$$
$$y\in \{yw^{n}| y\in Y(w), n\in\nN\} $$

{\footnotesize
\begin{table}[htbp] \caption{Condition {\bf I5}} 
\label{PcG}
\begin{center}
\renewcommand{\arraystretch}{1.4}
$\begin{array}{|c|c|c|} \hline
u\in U & x\in X_{u} & T_{x}C_{u} \mod \cH_{<0}\\ \hline

u_{1}^{4} & s_{1}s_{2}s_{1}s_{2}s_{3} & T_{xu} +v^{-a+3b}C_{u_{0}^{6}}\\\hline
u_{1}^{4} & s_{3}s_{1}s_{2}s_{1}s_{2}s_{3} & T_{xu} +v^{-a+3b}T_{s_{3}}C_{u_{0}^{6}}\\\hline
u_{1}^{4} & s_{2}s_{3}s_{1}s_{2}s_{1}s_{2}s_{3} & T_{xu} +v^{-a+3b}T_{s_{2}s_{3}}C_{u_{0}^{6}}\\\hline
u_{1}^{4} & s_{1}s_{2}s_{3}s_{1}s_{2}s_{1}s_{2}s_{3} & T_{xu} +v^{-a+3b}T_{s_{1}s_{2}s_{3}}C_{u_{0}^{6}}\\\hline
u_{1}^{4} & s_{2}s_{1}s_{2}s_{3}s_{1}s_{2}s_{1}s_{2}s_{3} & T_{xu} +v^{-a+3b}T_{s_{2}s_{1}s_{2}s_{3}}C_{u_{0}^{6}}\\\hline
u_{1}^{4} & yw & T_{xu} +v^{-a+3b}T_{y}C_{u_{0}^{6}}\\\hline

u_{1}^{6} & s_{1}s_{2}s_{1}s_{2}s_{3} & T_{xu} +v^{-a+3b}C_{u_{0}^{1}}\\
\hline

 \end{array}$
\end{center}
\end{table}}

\newpage


\item {\bf Case ${\bf r=2}$}

\begin{center}

\end{center}

\newpage

Let $w=s_{1}s_{2}s_{1}s_{2}s_{3}s_{1}s_{2}s_{1}s_{2}s_{3}$, $w_{1}=s_{1}s_{2}s_{3}s_{2}s_{1}s_{2}$ and $w_{2}=s_{2}s_{3}s_{2}s_{1}s_{2}s_{1}$. In the following table we have
$$y\in \{zw^{n}| z\in Y(w), n\in\nN\} $$
$$y_{1}\in \{yw_{1}^{n}| y\in Y(w_{1}), n\in\nN\} \text{ and } y_{2}\in \{yw_{2}^{n}| y\in Y(w_{2})\}$$

{\footnotesize
\begin{table}[htbp] \caption{Condition {\bf I5} for $r=2$} 
\label{PcG}
\begin{center}
\renewcommand{\arraystretch}{1.4}
$\begin{array}{|c|c|c|} \hline
u\in U & x\in X_{u} & T_{x}C_{u} \mod \cH_{<0}\\ \hline
u_{1}^{2} & ys_{1}s_{2}s_{1}s_{2}s_{3} & T_{xu} + T_{y}C_{u_{0}^{1}}\\\hline

u_{1}^{4} & s_{1}s_{2}s_{1}s_{2}s_{3} & T_{xu} +v C_{u_{0}^{6}}-C_{u_{0}^{9}}\\\hline
u_{1}^{4} & s_{3}s_{1}s_{2}s_{1}s_{2}s_{3} & T_{xu}+v T_{s_{3}}C_{u_{0}^{6}}-T_{s_{3}}C_{u_{0}^{9}}\\\hline
u_{1}^{4} & s_{2}s_{3}s_{1}s_{2}s_{1}s_{2}s_{3} & 
T_{xu} +v T_{s_{2}s_{3}}C_{u_{0}^{6}}-T_{s_{2}s_{3}}C_{u_{0}^{9}}\\\hline
u_{1}^{4} & s_{1}s_{2}s_{3}s_{1}s_{2}s_{1}s_{2}s_{3} & 
T_{xu} +v T_{s_{1}s_{2}s_{3}}C_{u_{0}^{6}}-T_{s_{1}s_{2}s_{3}}C_{u_{0}^{9}}\\\hline
u_{1}^{4} & s_{2}s_{1}s_{2}s_{3}s_{1}s_{2}s_{1}s_{2}s_{3} & 
T_{xu} +v T_{s_{2}s_{1}s_{2}s_{3}}C_{u_{0}^{6}}-T_{s_{2}s_{1}s_{2}s_{3}}C_{u_{0}^{9}}\\\hline
u_{1}^{4} & yw & 
T_{xu} +v T_{ys_{1}s_{2}s_{1}s_{2}s_{3}}C_{u_{0}^{6}}- T_{ys_{2}s_{1}s_{2}s_{1}s_{2}s_{3}}C_{u_{0}^{6}}-T_{ys_{1}s_{2}s_{1}s_{2}s_{3}}C_{u_{0}^{9}}\\
\hline

 \end{array}$
\end{center}
\end{table}}

\newpage


\item {\bf Case ${\bf 3/2<r<2}$}

\begin{center}

\end{center}

\vspace{2cm}

Let $w=s_{1}s_{2}s_{1}s_{2}s_{3}s_{1}s_{2}s_{1}s_{2}s_{3}$ and $w_{1}=s_{2}s_{3}s_{2}s_{1}s_{2}s_{1}$. In the following table we have
$$y\in \{zw^{n}| z\in Y(w), n\in\nN\} \text{ and } y_{1}\in \{yw_{1}^{n}| y\in Y(w_{1}), n\in\nN\} $$

\newpage

{\footnotesize
\begin{table}[htbp] \caption{Condition {\bf I5} $3/2<r<2$} 
\label{PcG}
\begin{textblock}{21}(0,4)
\renewcommand{\arraystretch}{1.3}
$\begin{array}{|c|c|c|} \hline
u\in U & x\in X_{u} & T_{x}C_{u} \mod \cH_{<0}\\ \hline
u_{1}^{2} & ys_{1}s_{2}s_{1}s_{2}s_{3} & T_{xu} +v^{-a+2b} T_{y}C_{u_{0}^{1}}\\\hline

u_{1}^{4} & s_{1}s_{2}s_{1}s_{2}s_{3} & T_{xu} +v^{-a+3b} C_{u_{0}^{6}}-v^{-a+2b}C_{u_{0}^{9}}\\\hline
u_{1}^{4} & s_{3}s_{1}s_{2}s_{1}s_{2}s_{3} & T_{xu} +v^{-a+3b} T_{s_{3}}C_{u_{0}^{6}}-v^{-a+2b}T_{s_{3}}C_{u_{0}^{9}}\\\hline
u_{1}^{4} & s_{2}s_{3}s_{1}s_{2}s_{1}s_{2}s_{3} & 
T_{xu} +v^{-a+3b} T_{s_{2}s_{3}}C_{u_{0}^{6}}-v^{-a+2b}T_{s_{2}s_{3}}C_{u_{0}^{9}}\\\hline
u_{1}^{4} & s_{1}s_{2}s_{3}s_{1}s_{2}s_{1}s_{2}s_{3} & 
T_{xu} +v^{-a+3b} T_{s_{1}s_{2}s_{3}}C_{u_{0}^{6}}-v^{-a+2b}T_{s_{1}s_{2}s_{3}}C_{u_{0}^{9}}\\\hline
u_{1}^{4} & s_{2}s_{1}s_{2}s_{3}s_{1}s_{2}s_{1}s_{2}s_{3} & 
T_{xu} +v^{-a+3b} T_{s_{2}s_{1}s_{2}s_{3}}C_{u_{0}^{6}}-v^{-a+2b}T_{s_{2}s_{1}s_{2}s_{3}}C_{u_{0}^{9}}\\\hline
u_{1}^{4} & yw& 
T_{xu} +v^{-a+3b} T_{ys_{1}s_{2}s_{1}s_{2}s_{3}}C_{u_{0}^{6}}-v^{-a+2b} T_{ys_{2}s_{1}s_{2}s_{1}s_{2}s_{3}}C_{u_{0}^{6}}-v^{-a+2b}T_{ys_{1}s_{2}s_{1}s_{2}s_{3}}C_{u_{0}^{9}}\\
\hline

 \end{array}$
 \end{textblock}
\end{table}}

\newpage


\item {\bf Case ${\bf r=3/2}$}

\begin{center}

\end{center}

\vspace{2cm}

Let $w=s_{1}s_{3}s_{2}s_{1}s_{2}s_{1}s_{3}s_{2}s_{1}s_{2}$, and $w_{1}=s_{2}s_{3}s_{2}s_{1}s_{2}s_{1}$. In the following table we have
$$y\in \{zw^{n}| z\in Y(w), n\in\nN\} \text{ and } y_{1}\in \{yw_{1}^{n}| y\in Y(w_{1}), n\in\nN\} $$

 $\ $\\
{\small
\begin{table}[htbp] \caption{Condition {\bf I5} for $r=3/2$} 
\label{PcG}
\begin{center}
\renewcommand{\arraystretch}{1.4}
$\begin{array}{|c|c|c|} \hline
u\in U & x\in X_{u} & T_{x}C_{u} \mod \cH_{<0}\\ \hline
u_{1}^{1} & s_{2}s_{1}s_{2}s_{1} & T_{xu} +v C_{u_{0}^{1}}\\\hline
u_{1}^{1} & s_{3}s_{2}s_{1}s_{2}s_{1} & T_{xu} +v T_{s_{3}}C_{u_{0}^{1}}\\\hline
u_{1}^{1} & y_{1}w_{1} & T_{xu} +v T_{ys_{2}s_{3}}C_{u_{0}^{1}}\\\hline
u_{1}^{1} & y_{1}w_{1}s_{3}s_{2}s_{1} & T_{xu}+vT_{ys_{3}s_{2}s_{1}s_{2}s_{3}}C_{u_{0}^{1}}\\\hline

u_{1}^{3} & s_{1}s_{2}s_{1} & T_{xu} +v C_{u_{0}^{6}}\\\hline
u_{1}^{3} & s_{2}s_{1}s_{2}s_{1} & T_{xu} +v^{3} C_{u_{0}^{6}}\\\hline
u_{1}^{3} & s_{3}s_{2}s_{1}s_{2}s_{1} & T_{xu} +v^{3} T_{s_{3}}C_{u_{0}^{6}}\\
\hline

 \end{array}$
\end{center}
\end{table}}

\newpage


\item {\bf Case ${\bf 1<r<3/2}$}
\begin{center}

\end{center}

\vspace{2cm}

Let $w=s_{1}s_{3}s_{2}s_{1}s_{2}s_{1}s_{3}s_{2}s_{1}s_{2}$, and $w_{1}=s_{2}s_{3}s_{2}s_{1}s_{2}s_{1}$. In the following table we have
$$y\in \{zw^{n}| z\in Y(w), n\in\nN\} \text{ and } y_{1}\in \{yw_{1}^{n}| y\in Y(w_{1}), n\in\nN\} $$

\newpage

{\footnotesize
\begin{table}[h!] \caption{Condition {\bf I5} for $1<r<3/2$} 
\label{PcG}
\begin{textblock}{21}(0,4)
\renewcommand{\arraystretch}{1.4}
$\begin{array}{|c|c|c|} \hline
u\in U & x\in X_{u} & T_{x}C_{u} \mod \cH_{<0}\\ \hline
u_{1}^{1} & s_{2}s_{1}s_{2}s_{1} & T_{xu} +v^{a+b} C_{u_{0}^{1}}\\\hline
u_{1}^{1} & s_{3}s_{2}s_{1}s_{2}s_{1} & T_{xu} +v^{a+b} T_{s_{3}}C_{u_{0}^{1}}\\\hline
u_{1}^{1} & y_{1}w_{1} & T_{xu} +v^{a+b} T_{y_{1}s_{2}s_{3}}C_{u_{0}^{1}}\\\hline
u_{1}^{1} & y_{1}w_{1}s_{3}s_{2}s_{1} & T_{xu}+v^{a+b}T_{y_{1}s_{3}s_{2}s_{1}s_{2}s_{3}}C_{u_{0}^{1}}\\\hline

u_{1}^{3} & s_{1}s_{2}s_{1} & T_{xu} +v^{a+b} C_{u_{0}^{6}}\\\hline
u_{1}^{3} & s_{2}s_{1}s_{2}s_{1} & T_{xu} +v^{a} C_{u_{0}^{6}}\\\hline
u_{1}^{3} & s_{3}s_{2}s_{1}s_{2}s_{1} & T_{xu} +v^{a} T_{s_{3}}C_{u_{0}^{6}}\\
\hline

 \end{array}$
 \end{textblock}
\end{table}}

$\ $\\

\newpage


\item {\bf Case ${\bf r=1}$}

\begin{center}

\end{center}

\newpage


\begin{center}


\end{center}

\vspace{2cm}

Let $w=s_{1}s_{3}s_{2}s_{1}s_{2}s_{1}s_{3}s_{2}s_{1}s_{2}$, and $w_{1}=s_{2}s_{3}s_{2}s_{1}s_{2}s_{1}$. In the following table we have
$$y\in \{zw^{n}| z\in Y(w), n\in\nN\} \text{ and } y_{1}\in \{yw_{1}^{n}| y\in Y(w_{1}), n\in\nN\} $$

\newpage

\vspace{2cm}
{\footnotesize
\begin{table}[t] \caption{Condition {\bf I5} for $r=1$} 
\label{PcG}
\begin{textblock}{21}(0,4)
\renewcommand{\arraystretch}{1.4}
$\begin{array}{|c|c|c|} \hline
u\in U & x\in X_{u} & T_{x}C_{u} \mod \cH_{<0}\\ \hline
u_{1}^{1} & s_{2}s_{1}s_{2}s_{1} & T_{xu} + C_{u_{0}^{1}}\\\hline
u_{1}^{1} & s_{3}s_{2}s_{1}s_{2}s_{1} & T_{xu} + T_{s_{3}}C_{u_{0}^{1}}\\\hline
u_{1}^{1} & y_{1}w_{1}& T_{xu} + T_{y_{1}s_{2}s_{3}}C_{u_{0}^{1}}\\\hline
u_{1}^{1} & y_{1}w_{1}s_{3}s_{2}s_{1} & T_{xu}+T_{y_{1}s_{3}s_{2}s_{1}s_{2}s_{3}}C_{u_{0}^{1}}\\\hline

u_{1}^{3} & s_{1}s_{2}s_{1} & T_{xu} + C_{u_{0}^{6}}\\\hline
u_{1}^{3} & s_{2}s_{1}s_{2}s_{1} & T_{xu} +v C_{u_{0}^{6}}\\\hline
u_{1}^{3} & s_{3}s_{2}s_{1}s_{2}s_{1} & T_{xu} +v T_{s_{3}}C_{u_{0}^{6}}\\
\hline

 \end{array}$
\end{textblock}
\end{table}}

$\ $\\

\newpage


\item {\bf Case ${\bf 0<r<1}$}

\begin{center}

\end{center}

\newpage

The green edges only hold for $1/2\leq r<1$.

\begin{center}
%
\end{center}

Let $w=s_{1}s_{3}s_{2}s_{1}s_{2}s_{1}s_{3}s_{2}s_{1}s_{2}$, and $w_{1}=s_{2}s_{3}s_{2}s_{1}s_{2}s_{1}$. In the following table we have
$$y\in \{zw^{n}| z\in Y(w), n\in\nN\} \text{ and } y_{1}\in \{yw_{1}^{n}| y\in Y(w_{1}), n\in\nN\} $$

\newpage


\newpage

\vspace{2cm}
{\small
\begin{table}[htbp] \caption{Condition {\bf I5} for $0<r<1$} 
\label{PcG}
\begin{center}
\renewcommand{\arraystretch}{1.4}
$
\begin{array}{|c|c|c|} \hline
u\in U & x\in X_{u} & T_{x}C_{u} \mod \cH_{<0}\\ \hline
u_{1}^{3} & s_{2}s_{1}s_{2}s_{1} & T_{xu} + v^{a}C_{u_{0}^{6}}\\\hline
u_{1}^{3} & s_{3}s_{2}s_{1}s_{2}s_{1} & T_{xu} + v^{a}T_{s_{3}}C_{u_{0}^{6}}\\\hline
u_{1}^{3} & y_{1}w_{1} & 
T_{xu} + v^{a}T_{y_{1}s_{2}s_{3}}C_{u_{0}^{6}}+T_{y_{1}s_{3}s_{2}s_{1}s_{2}s_{3}}C_{u_{0}^{6}}\\\hline

u_{1}^{3} & s_{3}s_{2}s_{1} & T_{xu} +T_{s_{2}s_{3}}C_{u_{0}^{6}}\\\hline
u_{1}^{3} & s_{1}s_{3}s_{2}s_{1} & T_{xu} +T_{s_{1}s_{2}s_{3}}C_{u_{0}^{6}}\\\hline
u_{1}^{3} & s_{2}s_{1}s_{3}s_{2}s_{1} & T_{xu} +T_{s_{2}s_{1}s_{2}s_{3}}C_{u_{0}^{6}}\\\hline
u_{1}^{3} & s_{1}s_{2}s_{1}s_{3}s_{2}s_{1} & 
T_{xu} +T_{s_{1}s_{2}s_{1}s_{2}s_{3}}C_{u_{0}^{6}}\\
\hline

 \end{array}$
\end{center}
\end{table}}

\end{longlist}
\newpage

\section{Partition of the affine Weyl group $\tB_{2}$ into cells}
\label{sec3}

\newcounter{b}
\setcounter{b}{16}

Let $(W,S)$ be the affine Weyl group of type $\tB_{2}$ with diagram and weight function given by
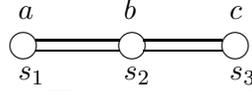
\begin{center}
%
\end{center}
where $a,b,c$ are positive integers. The partition of $W$ into cells only depend on the ratios $r_{1}=a/b$ and $r_{2}=c/b$. Note that by symmetry of diagram we may assume that $a\geq c$, that is $r_{1}\geq r_{2}$. We define the following open subsets of $\nR^{2}$:
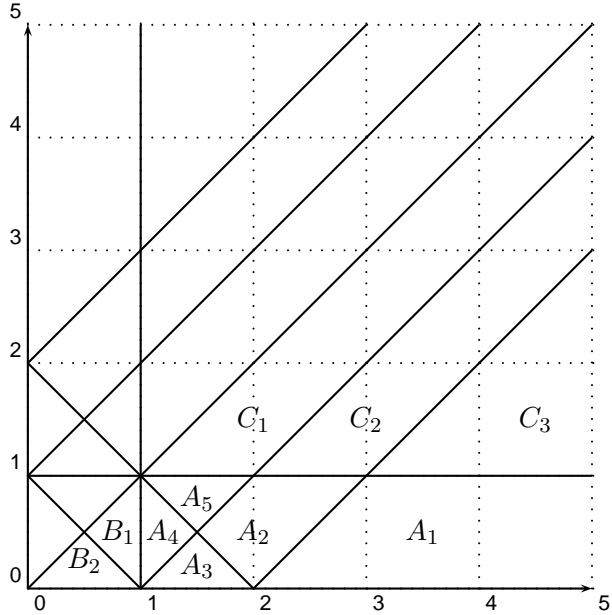
\begin{center}
\psset{xunit=1.5cm}
\psset{yunit=1.5cm}
%
\end{center}

\begin{Rem}
Let $V$ be an Euclidean space of dimension 3 with standard basis $v_{1},v_{2},v_{3}$ corresponding respectively to the conjugacy classes $\{s_{1}\}$, $\{s_{2}\}$ and $\{s_{3}\}$ in $S$.  In \cite{jeju4} we claimed that the essential hyperplanes of $W$ were 
$$\cH_{(1,0,0)},\cH_{(0,1,0)},\cH_{(0,0,1)},\cH_{(\eps,\eps,0)},\cH_{(0,\eps,\eps)},\cH_{(\eps,0,\eps)},\cH_{(\eps,\eps,\eps)},\cH_{(\eps,2\eps,\eps)}$$
where $\eps=\pm 1$. The picture above is just the projection of this hyperplane arrangement on the hyperplane $y=1$.
\end{Rem}

\newpage

\subsection{Partitions of $W$}
We describe the partitions $\cP_{\cLR,\fC}$ and $\cP_{\fC}$ of $\fC$ and the corresponding partition of $W$ according to \cite[Conjecture 3.3]{jeju4}. We keep the notation of  \cite[\S 5]{jeju4}.

\vspace{.5cm}

\begin{longlist}

\item Case $r_{1}=r_{2}, r_{2}>1$ \\

{\small
\begin{table}[htpb] \caption{Partition $\cP_{LR,\fC}$ and values of the $\ba$-function} 
\label{PlrcG}
\begin{center}
$\renewcommand{\arraystretch}{1.4}
\begin{array}{|l|c|} \hline
b_{8}=\{e\} & 0 \\
b_{7}=\{s_{2}\} & b\\
b_{6}=\{s_{1},s_{1}s_{2}, s_{2}s_{1},s_{2}s_{1}s_{2}\}&  a\\
b_{5}=\{s_{3}, s_{3}s_{2}, s_{2}s_{3}, s_{2}s_{3}s_{2}\}&  a\\
b_{4}=\{s_{3}s_{2}s_{3}\} & 2a-b \\
b_{3}=\{s_{1}s_{2}s_{1}\}& 2a-b\\ 
b_{2}=\{s_{1}s_{3}\}&  2a\\
b_{1}= \{w_{2,3}\} & 2a+2b\\ 
b_{0}= \{w_{1,2}\} & 2a+2b\\ \hline
\end{array}$
\end{center}
\end{table}
}

{\small
\begin{table}[h] \caption{Partition $\cP_{\fC}$} 
\label{PcG}
\begin{center}
\renewcommand{\arraystretch}{1.4}
$\begin{array}{|c|ccccccc|} \hline
r_{2}>1 & b_{0}\cup b_{1} & b_{2} & b_{3} & b_{4}  & b_{5}\cup b_{6} & b_{7} & b_{8}\\ \hline
 \end{array}$
\end{center}
\end{table}}

\psset{linewidth=.03mm}
\psset{unit=.5cm}

\begin{textblock}{10}(5.775,15.5)

\begin{center}

\end{center}

\end{textblock}

%


$\ $\\
\newpage

\item Case $r_{1}>1$ and $r_{2}>1$ \\


{\small
\begin{table}[h] \caption{Partition $\cP_{LR,\fC}$ and values of the $\ba$-function} 
\label{PcG}
\begin{center}
\renewcommand{\arraystretch}{1.4}
$\begin{array}{|l|c|} \hline
b_{8}=\{e\} & 0\\
b_{7}=\{s_{2}\} & b\\
b_{6}=\{s_{3},s_{2}s_{3},s_{3}s_{2},s_{2}s_{3}s_{2}\}   & c\\
b_{5}=\{s_{3}s_{2}s_{3}\} & 2c-b\\
b_{4}=\{s_{2}s_{3}s_{2}s_{3}\} & 2b+2c\\
b_{3}=\{s_{1},s_{2}s_{1},s_{1}s_{2},s_{2}s_{1}s_{2}\} & a\\
b_{2}=\{s_{1}s_{3}\}& a+c\\
b_{1}=\{s_{1}s_{2}s_{1}\} & 2a-b\\
b_{0}=\{w_{1,2}\}&  2a+2b\\ \hline
 \end{array}$
\end{center}
\end{table}}

{\small
\begin{table}[h] \caption{Partition $\cP_{\fC}$} 
\label{PcG}
\begin{center}
\renewcommand{\arraystretch}{1.4}
$\begin{array}{|c|ccccccccc|} \hline
(r_{1},r_{2})\in C_{3}  & b_{0} & b_{1} & b_{2} &\multicolumn{2}{r}{ b_{3}  \leftrightarrow b_{4}} & b_{5} & b_{6} & b_{7} & b_{8}\\ \hline
r_{1}-r_{2}-2=0  & b_{0} & b_{1} &\multicolumn{2}{r}{ b_{2} \cup b_{4}}  & b_{3} & b_{5} & b_{6} & b_{7} & b_{8}\\ \hline
(r_{1},r_{2})\in C_{2}  &b_{0} & \multicolumn{2}{r}{ b_{1}  \leftrightarrow b_{4}}& b_{2} & b_{3} & b_{5} & b_{6} & b_{7} & b_{8}\\ \hline 
r_{1}-r_{2}-1=0  & b_{0} & b_{4} &\multicolumn{2}{r}{ b_{1} \cup b_{2}}  & b_{3} & b_{5} & b_{6} & b_{7} & b_{8}\\ \hline
(r_{1},r_{2})\in C_{1}  &b_{0} &  b_{4}  & b_{2} & b_{1} & b_{3} & b_{5} & b_{6} & b_{7} & b_{8}\\ \hline 
\end{array}$
\end{center}
\end{table}}


\psset{linewidth=.13mm}

\begin{textblock}{10}(2,17)
\begin{center}

\end{center}

\end{textblock}


\vspace*{13.5cm}

\item Case $r_{1}>1$ and $r_{2}=1$\\
{\small
\begin{table}[h] \caption{Partition $\cP_{LR,\fC}$ and values of the $\ba$-function} 
\label{PlrcG}
\begin{center}
$\renewcommand{\arraystretch}{1.4}
\begin{array}{|l|c|} \hline
b_{6}=\{e\} & 0 \\
b_{5}=\{s_{2},s_{3},s_{3}s_{2},s_{2}s_{3},s_{3}s_{2}s_{3},s_{2}s_{3}s_{2}\} & b\\
b_{4}=\{s_{2}s_{3}s_{2}s_{3}\}& 4b\\ 
b_{3}=\{s_{1},s_{2}s_{1},s_{1}s_{2},s_{2}s_{1}s_{2}\}& a \\
b_{2}=\{s_{1}s_{3}\}&  a+b\\
b_{1}= \{s_{1}s_{2}s_{1}\} & 2a-b\\ 
b_{0}= \{w_{1,2}\} & 2a+2b\\ \hline
\end{array}$
\end{center}
\end{table}
}


\newpage

{\small
\begin{table}[h] \caption{Partition $\cP_{\fC}$} 
\label{PcG}
\begin{center}
\renewcommand{\arraystretch}{1.4}
$\begin{array}{|c|ccccccc|} \hline
r_{1}>3 & b_{0} & b_{1} & b_{2} &  \multicolumn{2}{r}{b_{3} \leftrightarrow b_{4}} & b_{5} & b_{6}\\ \hline

r_{1}=3 &b_{0} & b_{1} & \multicolumn{2}{r}{b_{2}\cup b_{4} }& b_{3} & b_{5} & b_{6}\\ \hline

3>r_{1}>2 &b_{0} & \multicolumn{2}{r}{b_{1} \leftrightarrow b_{4}} & b_{2}  & b_{3} & b_{5} & b_{6}\\ \hline

r_{1}=2 &b_{0} & b_{4} & \multicolumn{2}{r}{b_{1} \cup  b_{2}} & b_{3}  & b_{5} & b_{6}\\ \hline

2>r_{1}>1 &b_{0} & b_{4} & b_{2} &  b_{1} & b_{3}  & b_{5} & b_{6}\\ \hline

 \end{array}$
\end{center}
\end{table}}

\psset{linewidth=.13mm}

\begin{textblock}{10}(2,7.2)

\begin{center}

\end{center}

\end{textblock}


\newpage
\item Case $ r_{1}>1, r_{2}<1$ \\


{\small
\begin{table}[h!] \caption{Partition $\cP_{LR,\fC}$ and values of the $\ba$-function} 
\label{PlrcG}
\begin{center}
$\renewcommand{\arraystretch}{1.4}
\begin{array}{|l|c|} \hline
b_{8}=\{e\} & 0 \\
b_{7}=\{s_{3}\} & c\\
b_{6}=\{s_{2},s_{3}s_{2},s_{2}s_{3},s_{3}s_{2}s_{3}\}& b\\ 
b_{5}=\{s_{2}s_{3}s_{2}\} & 2b-c\\
b_{4}=\{s_{2}s_{3}s_{2}s_{3}\}& 2b+2c\\ 
b_{3}=\{s_{1},s_{2}s_{1},s_{1}s_{2},s_{2}s_{1}s_{2}\}& a \\
b_{2}=\{s_{1}s_{3}\}&  a+c\\
b_{1}= \{s_{1}s_{2}s_{1}\} & 2a-b\\ 
b_{0}= \{w_{1,2}\} & 2a+2b\\ \hline
\end{array}$
\end{center}
\end{table}
}


{\small
\begin{table}[h] \caption{Partition $\cP_{\fC}$} 
\label{PcG}
\begin{center}
\renewcommand{\arraystretch}{1.4}
$\begin{array}{|c|ccccccccc|} \hline
(r_{1},r_{2})\in A_{1}  & b_{0} & b_{1} & b_{2} &\multicolumn{2}{r}{ b_{3}  \leftrightarrow b_{4}} & b_{5} & b_{6} & b_{7} & b_{8}\\ \hline
r_{1}-r_{2}-2=0  & b_{0} & b_{1} &\multicolumn{2}{r}{ b_{2} \cup b_{4}}  & b_{3} & b_{5} & b_{6} & b_{7} & b_{8}\\ \hline
(r_{1},r_{2})\in A_{2}  &b_{0} & b_{1} & b_{4} & b_{2} & b_{3} & b_{5} & b_{6} & b_{7} & b_{8}\\ \hline 
r_{1}+r_{2}-2=0, r_{1}-r_{2}>1 &b_{0} & b_{1} & b_{4} & b_{2} &\multicolumn{2}{r}{ b_{3} \cup b_{5}} & b_{6} & b_{7} & b_{8}\\ \hline 
(r_{1},r_{2})\in A_{3} & b_{0} & \multicolumn{2}{r}{  b_{1} \leftrightarrow b_{4}}  & b_{2}  & b_{5} & b_{3} & b_{6} & b_{7} & b_{8} \\ \hline
r_{1}-r_{2}-1=0, r_{1}+r_{2}<2  & b_{0} &  b_{4} &  \multicolumn{2}{r}{b_{1} \cup b_{2}}  & b_{5} & b_{3} & b_{6} & b_{7} & b_{8} \\ \hline

(r_{1},r_{2})\in A_{4} & b_{0} & b_{4} & b_{2} & b_{1}  & b_{5} & b_{3} & b_{6} & b_{7} & b_{8} \\ \hline
r_{1}+r_{2}-2=0, r_{1}-r_{2}<1 & b_{0} & b_{4} & b_{2} & b_{1}  &\multicolumn{2}{r}{ b_{5} \cup b_{3}} & b_{6} & b_{7} & b_{8} \\ \hline
(r_{1},r_{2})\in A_{5} & b_{0} & b_{4} & b_{2} & b_{1} & b_{3} & b_{5} & b_{6} & b_{7} & b_{8} \\ \hline
r_{1}-r_{2}-1=0, r_{1}+r_{2}>2& b_{0} & b_{4} & \multicolumn{2}{r}{ b_{1} \cup b_{2}} & b_{3} & b_{5} & b_{6} & b_{7} & b_{8} \\ \hline
(r_{1},r_{2})=(3/2,1/2)  & b_{0} & b_{4} &\multicolumn{2}{r}{b_{1} \cup b_{2}} & \multicolumn{2}{r}{ b_{3} \cup b_{5}}& b_{6} & b_{7} & b_{8} \\ \hline
\end{array}$
\end{center}
\end{table}}

\psset{linewidth=.13mm}

\begin{textblock}{10}(2,19.5)

\begin{center}

\end{center}
\end{textblock}

 
$\ $\\ 
\newpage

\begin{textblock}{10}(5.775,2)
\begin{center}

\end{center}
\end{textblock}

\vspace*{7cm}

\item Case $ r_{1}=1, r_{2}<1$\\

{\small
\begin{table}[htpb] \caption{Partition $\cP_{LR,\fC}$ and values of the $\ba$-function} 
\label{PlrcG}
\begin{center}
$\renewcommand{\arraystretch}{1.4}
\begin{array}{|l|c|} \hline
b_{6}=\{e\} & 0 \\
b_{5}=\{s_{3}\} & c\\
b_{4}=\{s_{2},s_{1},s_{1}s_{2},s_{2}s_{1},s_{1}s_{2}s_{1},s_{3}s_{2},s_{2}s_{3},s_{3}s_{2}s_{3},s_{2}s_{1}s_{2}\} & a\\
b_{3}=\{s_{1}s_{3}\}&  a+c\\
b_{2}=\{s_{2}s_{3}s_{2}\} & 2a-c \\
b_{1}= \{w_{2,3}\} & 2a+2c\\ 
b_{0}= \{w_{1,2}\} & 4a\\ \hline
\end{array}$
\end{center}
\end{table}
}

{\small
\begin{table}[h] \caption{Partition $\cP_{\fC}$} 
\label{PcG}
\begin{center}
\renewcommand{\arraystretch}{1.4}
$\begin{array}{|c|ccccccc|} \hline
r_{1}=1, r_{2}<1 & b_{0} & b_{1} & \multicolumn{2}{r}{b_{2} \leftrightarrow b_{3}}  & b_{4} & b_{5} & b_{6}\\ \hline
 \end{array}$
\end{center}
\end{table}}


\begin{textblock}{10}(5.775,20)

\begin{center}

\end{center}

\end{textblock}


\newpage

\item Case $ r_{1}<1, r_{2}<1$

{\small
\begin{table}[h] \caption{Partition $\cP_{LR,\fC}$ and values of the $\ba$-function} 
\label{PlrcG}
\begin{center}
$\renewcommand{\arraystretch}{1.4}
\begin{array}{|l|c|} \hline
b_{8}=\{e\} & 0 \\
b_{7}=\{s_{3}\} & c\\
b_{6}=\{s_{1}\}& a\\ 
b_{5}=\{s_{1}s_{3}\}&  a+c\\
b_{4}=\{s_{2},s_{1}s_{2},s_{2}s_{1},s_{1}s_{2}s_{1},s_{3}s_{2},s_{2}s_{3},s_{3}s_{2}s_{3}\} & b\\
b_{3}=\{s_{2}s_{1}s_{2}\}& 2b-a\\ 
b_{2}=\{s_{2}s_{3}s_{2}\} & 2b-c \\
b_{1}= \{w_{2,3}\} & 2b+2c\\ 
b_{0}= \{w_{1,2}\} & 2a+2b\\ \hline
\end{array}$
\end{center}
\end{table}
}

\vspace{-.3cm}
{\small
\begin{table}[h] \caption{Partition $\cP_{\fC}$} 
\label{PcG}
\begin{center}
\renewcommand{\arraystretch}{1.4}
$\begin{array}{|c|ccccccccc|} \hline
(r_{1},r_{2})\in B_{2} & b_{0} & b_{1} & b_{2} & b_{3}  & b_{4} & b_{5} & b_{6} & b_{7} & b_{8}\\ \hline

r_{1}+r_{2}-1=0&b_{0} & b_{1} & b_{2} & b_{3} &  \multicolumn{2}{c}{b_{4} \cup b_{5}} & b_{6} & b_{7} & b_{8}\\ \hline

(r_{1},r_{2})\in B_{1}& b_{0} & b_{1} & b_{2} &\multicolumn{2}{l}{b_{3}  \leftrightarrow b_{5}} & b_{4} & b_{6} & b_{7} & b_{8}\\ \hline \end{array}$
\end{center}
\end{table}}

\psset{linewidth=.13mm}

\begin{textblock}{10}(2,14.2)

\begin{center}

\end{center}

\end{textblock}


\newpage

\item Case $r_{1}=r_{2}, r_{2}<1$\\

\vspace{-.4cm}

{\small
\begin{table}[htpb] \caption{Partition $\cP_{LR,\fC}$ and values of the $\ba$-function} 
\label{PlrcG}
\begin{center}
$\renewcommand{\arraystretch}{1.3}
\begin{array}{|l|c|} \hline
b_{8}=\{e\} & 0 \\
b_{7}=\{s_{3}\} & a\\
b_{6}=\{s_{1}\}&  a\\
b_{5}=\{s_{1}s_{3}\}&  2a\\
b_{4}=\{s_{2},s_{1}s_{2},s_{2}s_{1},s_{1}s_{2}s_{1},s_{3}s_{2},s_{2}s_{3},s_{3}s_{2}s_{3}\} & b \\
b_{3}=\{s_{2}s_{3}s_{2}\}& 2b-a\\ 
b_{2}=\{s_{2}s_{1}s_{2}\}& 2b-a\\ 
b_{1}= \{w_{2,3}\} & 2a+2b\\
b_{0}= \{w_{1,2}\} & 2a+2b\\ \hline
\end{array}$
\end{center}
\end{table}
}

\vspace{-.5cm}

{\small
\begin{table}[htbp] \caption{Partition $\cP_{\fC}$} 
\label{PcG}
\begin{center}
\renewcommand{\arraystretch}{1.3}
$\begin{array}{|c|ccccccc|} \hline
1>r_{2}>1/2 & (b_{0}\cup b_{1}) & (b_{2}\cup b_{3}) & b_{4} & b_{5}  & b_{6} & b_{7} & b_{8}\\ \hline

r_{2}=1/2&(b_{0}\cup b_{1}) & (b_{2}\cup b_{3}) & \multicolumn{2}{r}{b_{4}\cup b_{5} }& b_{6} & b_{7} & b_{8}\\ \hline

r_{2}<1/2&(b_{0}\cup b_{1}) & \multicolumn{2}{r}{(b_{2}\cup b_{3}) \overset{\frac{3}{2}}{\leftrightarrow} b_{5}} & b_{4}  & b_{6} & b_{7} & b_{8}\\ \hline
 \end{array}$
\end{center}
\end{table}}

\psset{linewidth=.13mm}

\begin{textblock}{10}(2,13.5)

\begin{center}

\end{center}

\end{textblock}

$\ $\\
\newpage

\item Case $ r_{1}=r_{2}=1$\\

{\small
\begin{table}[htbp] \caption{Partition $\cP_{LR,\fC}$ and values of the $\ba$-function} 
\label{PcG}
\begin{center}
\renewcommand{\arraystretch}{1.4}
$\begin{array}{|l|c|} \hline
b_{4}=\{e\} & 0\\
b_{3}=\fC-b_{0}\cup b_{1}\cup \{e\} & a\\
b_{2}=\{s_{1}s_{3}\}& 2a\\
b_{1}=\{w_{2,3}\}&  4a\\ 
b_{0}=\{w_{1,2}\}&  4a\\ \hline
 \end{array}$
\end{center}
\end{table}}

{\small
\begin{table}[h] \caption{Partition $\cP_{\fC}$} 
\label{PcG}
\begin{center}
\renewcommand{\arraystretch}{1.4}
$\begin{array}{|c|cccc|} \hline
r_{2}=r_{1}=1 & b_{0}\cup b_{1} & b_{2} & b_{3} & b_{4}\\ \hline
 \end{array}$
\end{center}
\end{table}}

\psset{linewidth=.13mm}

\begin{textblock}{10}(5.775,11)

\begin{center}

\end{center}
\end{textblock}


\end{longlist}

$\ $\\

\newpage 
\subsection{asymptotic cases}
We now consider the case where some parameters are equal to zero.\\

\begin{longlist}
\item Case $b=0, a,c>0$.  We have
$$W=(\nZ/2\nZ)\ltimes \big(\tA_{1}\times \tA_{1}\big)$$
where $\nZ/2\nZ$ is generated by $I=\{s_{2}\}$ and $\tA_{1}\times \tA_{1}$ is generated by $s_{1},s_{2}s_{3}s_{2},s_{3},s_{2}s_{1}s_{2}$. We know that the left (respectively two-sided) cells of $W$ are of the form $W_{I}.C$ (respectively $W_{I}.C.W_{I}$) where $C$ is a left (respectively two-sided) cell of $\tA_{1}\times \tA_{1}$ with respect to the weight function $\tL$ defined by 
$$\tL(s_{1})= \tL(s_{2}s_{1}s_{2})=a\text{ and } \tL(s_{3})=\tL(s_{2}s_{3}s_{2})=c.$$
We obtain the following partition of $W$. 

\psset{linewidth=.13mm}
\psset{unit=.5cm}

\begin{textblock}{10}(2,9)
\begin{center}

\end{center}

\end{textblock}


\vspace{7cm}

\item {\bf Case $c=0$}. \hspace{1cm} We have
$$W=(\nZ/2\nZ)\ltimes \tB_{2}$$
where $(\nZ/2\nZ)$ is generated by $s_{3}$ and $\tB_{2}$ is generated by $s_{2},s_{1},s_{3}s_{2}s_{3}$. We know that the left (respectively two-sided) cells of $W$ are of the form $W_{I}.C$ (respectively $W_{I}.C.W_{I}$) where $C$ is a left (respectively two-sided) cell of $\tB_{2}$ with respect to the weight function $\tL$ defined by 
$$\tL(s_{1})=a\text{ and } \tL(s_{2})=\tL(s_{3}s_{2}s_{3})=b.$$
We obtain the following partition of $W$.

\psset{linewidth=.13mm}

\begin{textblock}{10}(2,21)
\begin{center}

\end{center}

\end{textblock}

\end{longlist}

%
$\ $\\

\newpage

\section{Proof of \cite[Theorem 5.3]{jeju4} for $\tB_{2}$ }
\label{sec4}

As for $\tG_{2}$ (see Section 3), in each case we give the partial order $\preceq$ on $U$ (which is also the left order on the left cells) and a table for Condition {\bf I5}. When we have $T_{x}C_{u}\equiv T_{xu}$ we do not put it in the table. So when there is no table it means that $T_{x}C_{u}\equiv T_{xu}\mod \cH_{<0}$ for all $u\in U$ and $x\in X_{u}$.

\begin{longlist}

\item {\bf Case $(r_{1},r_{2})\in C_{3}$}
$\ $\\

\begin{center}

\end{center}
\end{textblock}


\newpage

\item {\bf Case $r_{1}-r_{2}-2=0$, $r_{2}>1$}
$\ $\\

\begin{center}

\end{center}


\newpage

\item {\bf Case $(r_{1},r_{2})\in C_{2}$}
$\ $\\

\begin{center}

\end{center}

\newpage

\item {\bf Case $r_{1}-r_{2}-1=0$, $r_{2}>1$}
$\ $\\

\begin{center}

\end{center}

\newpage

\item {\bf Case $(r_{1},r_{2})\in C_{1}$}
$\ $\\

\begin{center}

\end{center}

\vspace{.8cm}
Let $w=s_{1}s_{3}s_{2}s_{1}s_{3}s_{2}$ and $w_{1}=s_{1}s_{2}s_{3}s_{2}$. In the following table 
$$y\in \{zw^{n}|z\in Y(w), n\in \nN\} \text{ and } y_{1}\in \{zw^{n}|z\in Y(w_{1}), n\in \nN\} .$$
$\ $\\

{\small
\begin{table}[htbp] \caption{Condition {\bf I5} } 
\label{PcG}
\begin{center}
\renewcommand{\arraystretch}{1.4}
$\begin{array}{|c|c|c|} \hline
u\in U & x\in X_{u} & T_{x}C_{u} \mod \cH_{<0}\\ \hline
u_{2}^{1} & s_{1}s_{2} & T_{xu} +v^{a-c}\hspace{2pt} C_{u_{0}^{1}}\\\hline
u_{2}^{1} & ys_{1}s_{3}s_{2} & T_{xu} +v^{a-c}\hspace{2pt} T_{ys_{3}}C_{u_{0}^{1}}\\\hline
u_{4}^{3} & s_{2}s_{3}s_{2} & T_{xu}+v^{-a+2b}\hspace{2pt} C_{u_{1}^{3}}\\\hline
u_{4}^{3} & y_{1}s_{1}s_{2}s_{3}s_{2} & T_{xu} +v^{-a+2b}\hspace{2pt} T_{y_{1}s_{1}}C_{u_{1}^{3}}\\\hline

 \end{array}$
\end{center}
\end{table}}

\newpage

\item {\bf Case $ r_{1}=r_{2}, r_{2}>1$}

$\ $\\

\begin{center}

\end{center}

\vspace{1cm}

Let $w=s_{1}s_{3}s_{2}s_{1}s_{3}s_{2}$. In the following table
$$y\in \{zw^{n}| y\in Y(w), n\in\nN\}.$$

{\small
\begin{table}[h!] \caption{Condition {\bf I5}} 
\label{PcG}
\begin{center}
\renewcommand{\arraystretch}{1.4}
$\begin{array}{|c|c|c|} \hline
u\in U & x\in X_{u} & T_{x}C_{u} \mod \cH_{<0}\\ \hline
u_{2}^{1} & s_{1}s_{2} & T_{xu} +\hspace{2pt} C_{u_{0}^{1}}\\\hline
u_{2}^{1} & ys_{1}s_{3}s_{2} & T_{xu} +\hspace{2pt} T_{ys_{3}}C_{u_{0}^{1}}\\\hline
u_{2}^{3} & s_{3}s_{2} & T_{xu}+\hspace{2pt} C_{u_{0}^{6}}\\\hline
u_{2}^{3} & ys_{1}s_{3}s_{2} & T_{xu} +\hspace{2pt} T_{ys_{1}}C_{u_{0}^{6}}\\\hline

 \end{array}$
\end{center}
\end{table}}

$\ $\\


\newpage

\item{\bf Case $r_{1}>3,r_{2}=1$}
$\ $\\

\begin{center}

\end{center}

\newpage

\item{\bf Case $r_{1}=3,r_{2}=1$}
$\ $\\

\begin{center}

\end{center}

\newpage

\item{\bf Case $2<r_{1}<3,r_{2}=1$}
$\ $\\

\begin{center}

\end{center}

\newpage

\item{\bf Case $r_{1}=2,r_{2}=1$}
$\ $\\

\begin{center}

\end{center}

\vspace{2cm}

\newpage

Let $w=s_{1}s_{2}s_{3}s_{2}$. In the following table 
$$y\in \{zw^{n}|z\in Y(w), n\in \nN\}.$$
$\ $\\

{\small
\begin{table}[htbp] \caption{Condition {\bf I5}} 
\label{PcG}
\begin{center}
\renewcommand{\arraystretch}{1.4}
$\begin{array}{|c|c|c|} \hline
u\in U & x\in X_{u} & T_{x}C_{u} \mod \cH_{<0}\\ \hline
u_{3}^{3} & s_{2}s_{3}s_{2} & T_{xu}+ C_{u_{1}^{3}}\\\hline
u_{3}^{3} & ys_{1}s_{2}s_{3}s_{2} & T_{xu} + T_{ys_{1}}C_{u_{1}^{3}}\\\hline

 \end{array}$
\end{center}
\end{table}}

\newpage


\item{\bf Case $1<r_{1}<2,r_{2}=1$}
$\ $\\

\begin{center}

\end{center}

$\ $\\
\newpage

Let $w=s_{1}s_{3}s_{2}s_{1}s_{3}s_{2}$ and $w_{1}=s_{1}s_{2}s_{3}s_{2}$. In the following table 
$$y\in \{zw^{n}|z\in Y(w), n\in \nN\} \text{ and } y_{1}\in \{zw_{1}^{n}|z\in Y(w_{1}), n\in \nN\}$$

{\small
\begin{table}[htbp] \caption{Condition {\bf I5}} 
\label{PcG}
\begin{center}
\renewcommand{\arraystretch}{1.4}
$\begin{array}{|c|c|c|} \hline
u\in U & x\in X_{u} & T_{x}C_{u} \mod \cH_{<0}\\ \hline
u_{2}^{1} & s_{1}s_{2} & T_{xu} +v^{a-c}\hspace{2pt} C_{u_{0}^{1}}\\\hline
u_{2}^{1} & ys_{1}s_{3}s_{2} & T_{xu} +v^{a-c}\hspace{2pt} T_{ys_{3}}C_{u_{0}^{1}}\\\hline
u_{4}^{3} & s_{2}s_{3}s_{2} & T_{xu}+v^{-a+2c}\hspace{2pt} C_{u_{1}^{3}}\\\hline
u_{4}^{3} & y_{1}w_{1} & T_{xu} +v^{-a+2c}\hspace{2pt} T_{ys_{1}}C_{u_{1}^{3}}\\\hline

 \end{array}$
\end{center}
\end{table}}

\newpage


\item{\bf Case $(r_{1},r_{2})\in A_{1}$}
$\ $\\

\begin{center}

\end{center}

\newpage

\item{\bf Case $r_{1}-r_{2}-2=0$, $r_{2}<1$}
$\ $\\

\begin{center}

\end{center}

\newpage

\item{\bf Case $(r_{1},r_{2})\in A_{2}$}
$\ $\\

\begin{center}

\end{center}

\newpage

Let $w_{1}=s_{1}s_{2}s_{3}s_{2}$. In the following table 
$$y_{1}\in \{zw^{n}|z\in Y(w_{1}), n\in \nN\} $$
{\small
\begin{table}[htbp] \caption{Condition {\bf I5} } 
\label{PcG}
\begin{center}
\renewcommand{\arraystretch}{1.4}
$\begin{array}{|c|c|c|} \hline
u\in U & x\in X_{u} & T_{x}C_{u} \mod \cH_{<0}\\ \hline
u_{4}^{3} & s_{2}s_{3}s_{2} & T_{xu}+v^{-a+2b}\hspace{2pt} C_{u_{2}^{3}}\\\hline
u_{4}^{3} & y_{1}w_{1} & T_{xu} +v^{-a+2b}\hspace{2pt} T_{y_{1}s_{1}}C_{u_{2}^{3}}\\\hline

 \end{array}$
\end{center}
\end{table}}

$\ $\\

\newpage

\item{\bf Case $r_{1}+r_{2}-2=0, r_{1}-r_{2}>1$}
$\ $\\

\begin{center}

\end{center}

\newpage

\item{\bf Case $(r_{1},r_{2})\in A_{3}$}
$\ $\\

\begin{center}

\end{center}
\newpage

Let $w=s_{2}s_{3}s_{2}s_{1}$. In the following table 
$$y\in \{zw^{n}|z\in Y(w), n\in \nN\}.$$

{\small
\begin{table}[htbp] \caption{Condition {\bf I5} } 
\label{PcG}
\begin{center}
\renewcommand{\arraystretch}{1.4}
$\begin{array}{|c|c|c|} \hline
u\in U & x\in X_{u} & T_{x}C_{u} \mod \cH_{<0}\\ \hline
u_{4}^{2} & s_{2}s_{1} & T_{xu} -v^{c} C_{u_{0}^{4}}\\\hline
u_{4}^{2} & s_{3}s_{2}s_{1} & T_{xu} -v^{c} T_{s_{3}}C_{u_{0}^{4}}\\\hline
u_{4}^{2} & ys_{2}s_{3}s_{2}s_{1} & T_{xu}+T_{ys_{3}s_{2}s_{3}}C_{u_{0}^{4}} -v^{c} T_{ys_{2}s_{3}}C_{u_{0}^{4}} \\\hline
u_{4}^{4} & s_{2}s_{1} & T_{xu} -v^{c} C_{u_{0}^{1}}\\\hline
u_{4}^{4} & s_{3}s_{2}s_{1} & T_{xu} -v^{c} T_{s_{3}}C_{u_{0}^{1}}\\\hline
u_{4}^{4} & ys_{2}s_{3}s_{2}s_{1} & T_{xu}+T_{ys_{3}s_{2}s_{3}}C_{u_{0}^{1}} -v^{c} T_{ys_{2}s_{3}}C_{u_{0}^{1}} \\\hline

 \end{array}$
\end{center}
\end{table}}

\newpage

\item{\bf Case $r_{1}-r_{2}-1=0, r_{1}+r_{2}<2$}
$\ $\\

\begin{center}

\end{center}

\newpage

Let $w=s_{2}s_{3}s_{2}s_{1}$. In the following table 
$$y\in \{zw^{n}|z\in Y(w), n\in \nN\}.$$

{\small
\begin{table}[htbp] \caption{Condition {\bf I5} } 
\label{PcG}
\begin{center}
\renewcommand{\arraystretch}{1.4}
$\begin{array}{|c|c|c|} \hline
u\in U & x\in X_{u} & T_{x}C_{u} \mod \cH_{<0}\\ \hline
u_{3}^{2} & s_{2}s_{1} & T_{xu} -v^{c} C_{u_{0}^{4}}\\\hline
u_{3}^{2} & s_{3}s_{2}s_{1} & T_{xu} -v^{c} T_{s_{3}}C_{u_{0}^{4}}\\\hline
u_{3}^{2} & ys_{2}s_{3}s_{2}s_{1} & T_{xu}+T_{ys_{3}s_{2}s_{3}}C_{u_{0}^{4}} -v^{c} T_{ys_{2}s_{3}}C_{u_{0}^{4}} \\\hline
u_{3}^{4} & s_{2}s_{1} & T_{xu} -v^{c} C_{u_{0}^{1}}\\\hline
u_{3}^{4} & s_{3}s_{2}s_{1} & T_{xu} -v^{c} T_{s_{3}}C_{u_{0}^{1}}\\\hline
u_{3}^{4} & ys_{2}s_{3}s_{2}s_{1} & T_{xu}+T_{ys_{3}s_{2}s_{3}}C_{u_{0}^{1}} -v^{c} T_{ys_{2}s_{3}}C_{u_{0}^{1}} \\\hline

 \end{array}$
\end{center}
\end{table}}

\newpage

\item{\bf Case $(r_{1},r_{2})\in A_{4}$}
$\ $\\

\begin{center}


\end{center}

$\ $\\

In the following diagram, the green edges only hold for $a\geq 2c$.

\begin{center}
\psset{unit=.62cm}
%
\end{center}

\newpage

Let $w=s_{1}s_{3}s_{2}s_{1}s_{3}s_{2}$ and $w_{1}=s_{2}s_{3}s_{2}$. In the following table 
$$y\in \{zw^{n}|z\in Y(w), n\in \nN\} \text{ and } y_{1}\in \{zw^{n}|z\in Y(w_{1}), n\in \nN\}$$

{\footnotesize
\begin{table}[htbp] \caption{Condition {\bf I5} } 
\label{PcG}
\begin{center}
\renewcommand{\arraystretch}{1.4}
$\begin{array}{|c|c|c|} \hline
u\in U & x\in X_{u} & T_{x}C_{u} \mod \cH_{<0}\\ \hline
u_{2}^{1} & s_{1}s_{2} & T_{xu} +v^{a-c} C_{u_{0}^{1}}\\\hline
u_{2}^{1} & ys_{1}s_{3}s_{2} & T_{xu}+v^{a-c} T_{ys_{3}}C_{u_{0}^{1}} \\\hline
u_{4}^{2} & s_{2}s_{1} & T_{xu} -v^{c} C_{u_{0}^{4}}\\\hline
u_{4}^{2} & s_{3}s_{2}s_{1} & T_{xu} -v^{c} T_{s_{3}}C_{u_{0}^{4}}\\\hline
u_{4}^{2} & y_{1}w_{1} & T_{xu}+T_{y_{1}s_{3}s_{2}s_{3}}C_{u_{0}^{4}} -v^{c} T_{y_{1}s_{2}s_{3}}C_{u_{0}^{4}} \\\hline
u_{4}^{4} & s_{2}s_{1} & T_{xu} -v^{c} C_{u_{0}^{1}}\\\hline
u_{4}^{4} & s_{3}s_{2}s_{1} & T_{xu} -v^{c} T_{s_{3}}C_{u_{0}^{1}}\\\hline
u_{4}^{4} & y_{1}w_{1}& T_{xu}+T_{y_{1}s_{3}s_{2}s_{3}}C_{u_{0}^{1}} -v^{c} T_{y_{1}s_{2}s_{3}}C_{u_{0}^{1}} \\\hline

 \end{array}$
\end{center}
\end{table}}

\newpage

\item{\bf Case $r_{1}+r_{2}-2=0, r_{1}-r_{2}<1$}
$\ $\\

\begin{center}


\end{center}

$\ $\\

In the following diagram, the green edges only hold for $a\geq 2c$.

\begin{center}
\psset{unit=.65cm}
%
\end{center}

\newpage
Let $w=s_{1}s_{3}s_{2}s_{1}s_{3}s_{2}$. In the following table 
$$y\in \{zw^{n}|z\in Y(w), n\in \nN\}.$$

{\small
\begin{table}[htbp] \caption{Condition {\bf I5} } 
\label{PcG}
\begin{center}
\renewcommand{\arraystretch}{1.4}
$\begin{array}{|c|c|c|} \hline
u\in U & x\in X_{u} & T_{x}C_{u} \mod \cH_{<0}\\ \hline
u_{2}^{1} & s_{1}s_{2} & T_{xu} +v^{a-c}\hspace{2pt} C_{u_{0}^{1}}\\\hline
u_{2}^{1} & ys_{1}s_{3}s_{2} & T_{xu} +v^{a-c}\hspace{2pt} T_{ys_{3}}C_{u_{0}^{1}}\\\hline

 \end{array}$
\end{center}
\end{table}}

\newpage

\item{\bf Case $(r_{1},r_{2})\in A_{5}$}
$\ $\\

\begin{center}


\end{center}

$\ $\\

In the following diagram, the green edges only hold for $a\geq 2c$.

\begin{center}
\psset{unit=.65cm}
%
\end{center}

\newpage 

Let $w=s_{1}s_{3}s_{2}s_{1}s_{3}s_{2}$ and $w_{1}=s_{1}s_{2}s_{3}s_{2}$. In the following table 
$$y\in \{zw^{n}|z\in Y(w), n\in \nN\} \text{ and } y_{1}\in \{zw_{1}^{n}|z\in Y(w_{1}), n\in \nN\}$$

{\small
\begin{table}[htbp] \caption{Condition {\bf I5}} 
\label{PcG}
\begin{center}
\renewcommand{\arraystretch}{1.4}
$\begin{array}{|c|c|c|} \hline
u\in U & x\in X_{u} & T_{x}C_{u} \mod \cH_{<0}\\ \hline
u_{2}^{1} & s_{1}s_{2} & T_{xu} +v^{a-c}\hspace{2pt} C_{u_{0}^{1}}\\\hline
u_{2}^{1} & ys_{1}s_{3}s_{2} & T_{xu} +v^{a-c}\hspace{2pt} T_{ys_{3}}C_{u_{0}^{1}}\\\hline
u_{4}^{3} & s_{2}s_{3}s_{2} & T_{xu}+v^{-a+2b}\hspace{2pt} C_{u_{1}^{3}}\\\hline
u_{4}^{3} & y_{1}w_{1} & T_{xu} +v^{-a+2b}\hspace{2pt} T_{y_{1}s_{1}}C_{u_{1}^{3}}\\\hline

 \end{array}$
\end{center}
\end{table}}

\newpage

\item{\bf Case $r_{1}-r_{2}-1=0, r_{1}+r_{2}>2$,  $r_{2}<1$}
$\ $\\

\begin{center}

\end{center}

\newpage 

Let $w_{1}=s_{1}s_{2}s_{3}s_{2}$. In the following table 
$$ y_{1}\in \{zw_{1}^{n}|z\in Y(w_{1}), n\in \nN\}$$

{\small
\begin{table}[htbp] \caption{Condition {\bf I5}} 
\label{PcG}
\begin{center}
\renewcommand{\arraystretch}{1.4}
$\begin{array}{|c|c|c|} \hline
u\in U & x\in X_{u} & T_{x}C_{u} \mod \cH_{<0}\\ \hline
u_{3}^{3} & s_{2}s_{3}s_{2} & T_{xu} +v^{b-c}\hspace{2pt} C_{u_{1}^{3}}\\\hline
u_{3}^{3} & y_{1}s_{1}s_{2}s_{3}s_{2} & T_{xu} +v^{b-c}\hspace{2pt} T_{y_{1}s_{1}}C_{u_{1}^{3}}\\\hline

 \end{array}$
\end{center}
\end{table}}


\newpage

\item{\bf Case $(r_{1},r_{2})=(3/2,1/2)$}
$\ $\\

\begin{center}

\end{center}

\newpage

\item{\bf Case $r_{1}=1, r_{2}<1$}
$\ $\\

\begin{center}


\end{center}

\vspace{2cm} 

In the following diagram, the green edges only hold for $a\geq 2c$.

\vspace{1cm}

\begin{center}
\psset{unit=.65cm}
%
\end{center}

\newpage

Let $w=s_{1}s_{3}s_{2}s_{1}s_{3}s_{2}$ and $w_{1}=s_{2}s_{3}s_{2}s_{1}$. In the following table 
$$y\in \{zw^{n}|z\in Y(w), n\in \nN\} \text{ and } y_{1}\in \{zw^{n}|z\in Y(w_{1}), n\in \nN\}$$

{\small
\begin{table}[htbp] \caption{Condition {\bf I5}} 
\label{PcG}
\begin{center}
\renewcommand{\arraystretch}{1.4}
$\begin{array}{|c|c|c|} \hline
u\in U & x\in X_{u} & T_{x}C_{u} \mod \cH_{<0}\\ \hline
u_{2}^{2} & s_{2}s_{1} & T_{xu} -v^{c} C_{u_{0}^{4}}\\\hline
u_{2}^{2} & s_{3}s_{2}s_{1} & T_{xu} -v^{c} T_{s_{3}}C_{u_{0}^{4}}\\\hline
u_{2}^{2} & y_{1}w_{1} & T_{xu}+T_{y_{1}s_{3}s_{2}s_{3}}C_{u_{0}^{4}} -v^{c} T_{y_{1}s_{2}s_{3}}C_{u_{0}^{4}} \\\hline
u_{2}^{4} & s_{2}s_{1} & T_{xu} -v^{c} C_{u_{0}^{1}}\\\hline
u_{2}^{4} & s_{3}s_{2}s_{1} & T_{xu} -v^{c} T_{s_{3}}C_{u_{0}^{1}}\\\hline
u_{2}^{4} & ys_{2}s_{3}s_{2}s_{1} & T_{xu}+T_{ys_{3}s_{2}s_{3}}C_{u_{0}^{1}} -v^{c} T_{ys_{2}s_{3}}C_{u_{0}^{1}} \\\hline
u_{3}^{1} & s_{1}s_{2} & T_{xu} +v^{a-c} C_{u_{0}^{1}}\\\hline
u_{3}^{1} & ys_{1}s_{3}s_{2} & T_{xu} +v^{a-c}T_{ys_{3}}C_{u_{0}^{1}}\\\hline

 \end{array}$
\end{center}
\end{table}}


\newpage

\item{Case $(r_{1},r_{2})\in B_{1}$}
$\ $\\

\begin{center}


\end{center}

\vspace{2cm} 

In the following diagram, the green edges only hold for $a\geq 2c$.

\begin{center}
\psset{unit=.65cm}
%
\end{center}

\newpage

Let $w=s_{1}s_{3}s_{2}s_{1}s_{3}s_{2}$ and $w_{1}=s_{2}s_{3}s_{2}s_{1}$. In the following table 
$$y\in \{zw^{n}|z\in Y(w), n\in \nN\} \text{ and } y_{1}\in \{zw^{n}|z\in Y(w_{1}), n\in \nN\} .$$

{\small
\begin{table}[htbp] \caption{Condition {\bf I5} } 
\label{PcG}
\begin{center}
\renewcommand{\arraystretch}{1.4}
$\begin{array}{|c|c|c|} \hline
u\in U & x\in X_{u} & T_{x}C_{u} \mod \cH_{<0}\\ \hline
u_{2}^{2} & s_{2}s_{1} & T_{xu} -v^{c} C_{u_{0}^{4}}\\\hline
u_{2}^{2} & s_{3}s_{2}s_{1} & T_{xu} -v^{c} T_{s_{3}}C_{u_{0}^{4}}\\\hline
u_{2}^{2} & y_{1}w_{1} & T_{xu}+T_{y_{1}s_{3}s_{2}s_{3}}C_{u_{0}^{4}} -v^{c} T_{y_{1}s_{2}s_{3}}C_{u_{0}^{4}} \\\hline
u_{2}^{4} & s_{2}s_{1} & T_{xu} -v^{c} C_{u_{0}^{1}}\\\hline
u_{2}^{4} & s_{3}s_{2}s_{1} & T_{xu} -v^{c} T_{s_{3}}C_{u_{0}^{1}}\\\hline
u_{2}^{4} & y_{1}w_{1} & T_{xu}+T_{y_{1}s_{3}s_{2}s_{3}}C_{u_{0}^{1}} -v^{c} T_{y_{1}s_{2}s_{3}}C_{u_{0}^{1}} \\\hline
u_{3}^{1} & s_{1}s_{2} & T_{xu} +v^{a-c} C_{u_{0}^{1}}\\\hline
u_{3}^{1} & s_{1}s_{3}s_{2} & T_{xu} +v^{a-c}T_{s_{3}}C_{u_{0}^{1}}\\\hline
u_{3}^{1} & s_{2}s_{1}s_{3}s_{2} & T_{xu} +v^{a-c}T_{s_{2}s_{3}}C_{u_{0}^{1}}\\\hline
u_{3}^{1} & s_{3}s_{2}s_{1}s_{3}s_{2} & T_{xu} +v^{a-c}T_{s_{3}s_{2}s_{3}}C_{u_{0}^{1}}\\\hline
u_{3}^{1} & s_{1}s_{2}s_{1}s_{3}s_{2} & T_{xu} +v^{a-c}T_{s_{1}s_{2}s_{3}}C_{u_{0}^{1}}\\\hline
u_{3}^{1} & ys_{1}s_{3}s_{2}s_{1}s_{3}s_{2} & T_{xu} +v^{a-c}T_{ys_{1}s_{3}s_{2}s_{3}}C_{u_{0}^{1}}\\\hline

 \end{array}$
\end{center}
\end{table}}


\newpage

\item{\bf Case $r_{1}+r_{2}-1=0$}
$\ $\\

\begin{center}


\end{center}

$\ $\\

In the following diagram, the green edges only hold for $a\geq 2c$.

\begin{center}
\psset{unit=.65cm}
%
\end{center}

\newpage

Let $w_{1}=s_{2}s_{3}s_{2}s_{1}$. In the following table 
$$y_{1}\in \{zw^{n}|z\in Y(w_{1}), n\in \nN\}.$$

{\small
\begin{table}[htbp] \caption{Condition {\bf I5} } 
\label{PcG}
\begin{center}
\renewcommand{\arraystretch}{1.4}
$\begin{array}{|c|c|c|} \hline
u\in U & x\in X_{u} & T_{x}C_{u} \mod \cH_{<0}\\ \hline
u_{2}^{2} & s_{2}s_{1} & T_{xu} -v^{c} C_{u_{0}^{4}}\\\hline
u_{2}^{2} & s_{3}s_{2}s_{1} & T_{xu} -v^{c} T_{s_{3}}C_{u_{0}^{4}}\\\hline
u_{2}^{2} & ys_{2}s_{3}s_{2}s_{1} & T_{xu}+T_{ys_{3}s_{2}s_{3}}C_{u_{0}^{4}} -v^{c} T_{ys_{2}s_{3}}C_{u_{0}^{4}} \\\hline
u_{2}^{4} & s_{2}s_{1} & T_{xu} -v^{c} C_{u_{0}^{1}}\\\hline
u_{2}^{4} & s_{3}s_{2}s_{1} & T_{xu} -v^{c} T_{s_{3}}C_{u_{0}^{1}}\\\hline
u_{2}^{4} & ys_{2}s_{3}s_{2}s_{1} & T_{xu}+T_{ys_{3}s_{2}s_{3}}C_{u_{0}^{1}} -v^{c} T_{ys_{2}s_{3}}C_{u_{0}^{1}} \\\hline

 \end{array}$
\end{center}
\end{table}}

\newpage

\item{\bf Case $(r_{1},r_{2})\in B_{2}$}
$\ $\\

\begin{center}


\end{center}

$\ $\\

In the following diagram, the green edges only hold for $a\geq 2c$.

\begin{center}
\psset{unit=.65cm}
%
\end{center}

\newpage

Let $w_{1}=s_{2}s_{3}s_{2}s_{1}$. In the following table 
$$y_{1}\in \{zw_{1}^{n}|z\in Y(w_{1}), n\in \nN\}.$$

{\small
\begin{table}[htbp] \caption{Condition {\bf I5}  } 
\label{PcG}
\begin{center}
\renewcommand{\arraystretch}{1.4}
$\begin{array}{|c|c|c|} \hline
u\in U & x\in X_{u} & T_{x}C_{u} \mod \cH_{<0}\\ \hline
u_{2}^{2} & s_{2}s_{1} & T_{xu} -v^{c} C_{u_{0}^{4}}\\\hline
u_{2}^{2} & s_{3}s_{2}s_{1} & T_{xu} -v^{c} T_{s_{3}}C_{u_{0}^{4}}\\\hline
u_{2}^{2} & ys_{2}s_{3}s_{2}s_{1} & T_{xu}+T_{ys_{3}s_{2}s_{3}}C_{u_{0}^{4}} -v^{c} T_{ys_{2}s_{3}}C_{u_{0}^{4}} \\\hline
u_{2}^{4} & s_{2}s_{1} & T_{xu} -v^{c} C_{u_{0}^{1}}\\\hline
u_{2}^{4} & s_{3}s_{2}s_{1} & T_{xu} -v^{c} T_{s_{3}}C_{u_{0}^{1}}\\\hline
u_{2}^{4} & ys_{2}s_{3}s_{2}s_{1} & T_{xu}+T_{ys_{3}s_{2}s_{3}}C_{u_{0}^{1}} -v^{c} T_{ys_{2}s_{3}}C_{u_{0}^{1}} \\\hline

 \end{array}$
\end{center}
\end{table}}

\newpage

\item{\bf Case $r_{1}=r_{2}, 1/2<r_{2}<1$}
$\ $\\

\begin{center}

\end{center}

\newpage

Let $w=s_{1}s_{3}s_{2}s_{1}s_{3}s_{2}$. In the following table 
$$y\in \{zw^{n}|z\in Y(w), n\in \nN\}.$$

{\small
\begin{table}[htbp] \caption{Condition {\bf I5}} 
\label{PcG}
\begin{center}
\renewcommand{\arraystretch}{1.4}
$\begin{array}{|c|c|c|} \hline
u\in U & x\in X_{u} & T_{x}C_{u} \mod \cH_{<0}\\ \hline
u_{2}^{1} & s_{1}s_{2} & T_{xu} -+C_{u_{0}^{1}}\\\hline
u_{2}^{1} & ys_{1}s_{3}s_{2} & T_{xu} +T_{ys_{3}}C_{u_{0}^{1}}\\\hline
u_{2}^{3} & s_{3}s_{2} & T_{xu} +C_{u_{0}^{6}}\\\hline
u_{2}^{3} & ys_{1}s_{3}s_{2} & T_{xu} +T_{ys_{1}}C_{u_{0}^{6}}\\\hline

 \end{array}$
\end{center}
\end{table}}


\newpage

\item{\bf Case $r_{1}=r_{2}, r_{2}=1/2$}
$\ $\\

\begin{center}

\end{center}

$\ $\\
\newpage

\item{\bf Case $r_{1}=r_{2}, 0<r_{2}<1/2$}
$\ $\\

\begin{center}

\end{center}

\newpage

\item{\bf Case $r_{1}=r_{2}=1$}
$\ $\\

\begin{center}

\end{center}
$\ $\\

\newpage

Let $w=s_{1}s_{3}s_{2}s_{1}s_{3}s_{2}$. In the following table 
$$y\in \{zw^{n}|z\in Y(w), n\in \nN\}.$$

{\small
\begin{table}[htbp] \caption{Condition {\bf I5} } 
\label{PcG}
\begin{center}
\renewcommand{\arraystretch}{1.4}
$\begin{array}{|c|c|c|} \hline
u\in U & x\in X_{u} & T_{x}C_{u} \mod \cH_{<0}\\ \hline
u_{2}^{1} & s_{1}s_{2} & T_{xu} +\hspace{2pt} C_{u_{0}^{1}}\\\hline
u_{2}^{1} & ys_{1}s_{3}s_{2} & T_{xu} +\hspace{2pt} T_{ys_{3}}C_{u_{0}^{1}}\\\hline
u_{2}^{3} & s_{3}s_{2} & T_{xu}+\hspace{2pt} C_{u_{0}^{6}}\\\hline
u_{2}^{3} & ys_{1}s_{3}s_{2} & T_{xu} +\hspace{2pt} T_{ys_{1}}C_{u_{0}^{6}}\\\hline

 \end{array}$
\end{center}
\end{table}}

\end{longlist}

\end{document}